\documentclass[11pt,a4paper]{amsart}

\usepackage[pdftex, 
	a4paper,
	unicode,
	linktocpage,
	urlbordercolor={1 0.5 0.125},
	filebordercolor={1 0.5 0.125},       
	linkbordercolor={0.25 1 0.25},
	citebordercolor={0.25 1 0.25},
	pdfusetitle,
	pagebackref=false,
	pdfstartview=FitH,
	pdfborderstyle={/S/U/W 1}, 
	breaklinks]{hyperref}

\usepackage[T1]{fontenc}
\usepackage{lmodern}
\usepackage[utf8]{inputenc}
\usepackage[english,german,french]{babel}

\usepackage{graphicx}
\usepackage{amssymb}

\usepackage{verbatim} 

\newtheorem*{TheoremA}{Théorème A}
\newtheorem*{TheoremB}{Théorème B}
\newtheorem*{TheoremC}{Théorème C}
\newtheorem{theorem}{Théorème}[section]
\newtheorem{proposition}[theorem]{Proposition}
\newtheorem{lemma}[theorem]{Lemme}
\newtheorem{corollary}[theorem]{Corollaire}

\theoremstyle{definition}
\newtheorem*{definition}{Définition}
\newtheorem{question}{Question}

\theoremstyle{remark}
\newtheorem{remark}[theorem]{Remarque}

\newcommand{\ie}{i\textup.e\textup.{}}

\newcommand{\Int}[2]{\mathfrak{Int}_{#1}(#2)}
\newcommand{\id}{\operatorname{id}}
\newcommand{\fix}{\operatorname{Fix}}
\newcommand{\supp}{\operatorname{Supp}}

\newcommand{\Th}{\operatorname{Th}}
\newcommand{\Mod}{\operatorname{Mod}}
\newcommand{\arity}{\operatorname{ar}}
\newcommand{\V}{\mathcal V}
\newcommand{\T}{\mathcal T}
\newcommand{\F}{\mathcal F}

\title[Arithm\'etique dans des groupes
de permutations d'un intervalle]
{Interpr\'etation de l'Arithm\'etique\\
dans certains groupes de permutations\\
affines par morceaux d'un intervalle}

\author[T.~Alt{\i}nel et A.~Muranov]{Tuna Alt{\i}nel et Alexey Muranov}

\date{\today}


\address{Universit\'e de Lyon ;
Universit\'e Lyon 1 ;
Institut Camille Jordan CNRS UMR 5208 ;
43 boulevard du 11 novembre 1918,
F--69622 Villeurbanne Cedex, France}


\thanks{Le second auteur a \'et\'e soutenu d'abord par
une bourse Chateaubriand et apr\`es par
\emph{FP6 Marie Curie Research Training Network
in Model Theory and its Applications\/}
financ\'e par la commission europ\'eenne
sous le contrat num\'ero \texttt{MRTN-CT-2004-512234} (MODNET)}

\begin{document}


\maketitle

{\selectlanguage{english}
\begin{abstract}
The Arithmetic is interpreted in all the groups of Richard Thompson
and Graham Higman, as well as in other groups of piecewise affine
permutations of an interval which generalize the groups of Thompson
and Higman.
In particular, the elementary theories of all these groups are
undecidable.
Moreover, Thompson's group $F$ and some of its generalizations
interpret the Arithmetic without parameters.
\end{abstract}
}


\begin{abstract}
L'Arithmétique est interprétée dans tous les groupes de Richard Thompson
et de Graham Higman, aussi bien que dans d'autres groupes
des permutations affines par morceaux d'un intervalle qui
généralisent les groupes de Thompson et de Higman.
En particulier, les théories élémentaires de tous ces groupes sont
indécidables.
De plus, le groupe $F$ de Thompson et certaines de ses généralisations
interprètent l'Arithmétique sans paramètres.
\end{abstract}


\tableofcontents



\section{Introduction}
\label{section:introduction}

Valery Bardakov et Vladimir Tolstykh
\cite{BardakovTolstykh:2007:iaTgF} ont récemment montré que
le groupe $F$ de Richard Thompson interprète l'Arithmétique.
En d'autres termes, $F$ interprète la structure $(\mathbb N,+,\times)$
par des formules du premier ordre avec paramètres.
Dans ce travail nous généralisons ce résultat dans deux directions.
D'un côté, dans les sections \ref{section:twolemmas} et
\ref{section:copiesZwrZ}
nous généralisons l'approche de Bardakov et Tolstykh à
la fabrication de l'Arithmétique à partir de $F$ et montrons
qu'elle marche pour tous les groupes définis par
Melanie Stein dans \cite{Stein:1992:gplh},
dont tous les trois groupes de Thompson et tous les groupes de
Graham Higman \cite{Higman:1974:fpisg}.
D'un autre côté, dans la section \ref{section:interpretationarithmetic}
nous démontrons que le groupe $F$ et certaines de ses
généralisations interprètent l'Arithmétique \emph{sans paramètres}.
(La différence entre \emph{avec\/} et \emph{sans\/} paramètres
sera expliquée dans la section~\ref{subsection:generalities.modeltheory}.)


La théorie élémentaire de l'Arithmétique $(\mathbb N,+,\times)$
est notoire pour sa complexité depuis
les théorèmes d'incomplétude de Gödel
\cite{Godel:2006:fuSPMvS1-ger}.
L'une des significations de l'interprétabilité de l'Arithmétique dans
une structure de signature finie est ce qu'elle entraîne
l'\emph{indécidabilité héréditaire\/}
de la théorie élémentaire de cette structure.
(Une théorie de signature finie est dite
\emph{héréditairement indécidable\/}
si toute sous-théorie de la même signature est indécidable,
voir \cite[\S3]{Tarski:1971:gmpu}.)
Grâce à des travaux d'Andrzej Mostowski, de Raphael Robinson,
et d'Alfred Tarski
\cite{MostowskiRobTar:1971:ueua,MostowskiTarski:1949:uaitr,
Tarski:1971:gmpu},
il est bien connu que la théorie élémentaire de l'Arithmétique est
héréditairement indécidable.
Il est également bien connu aux spécialistes
que si une structure $N$ de signature finie interprète avec paramètres
une autre structure $M$ de signature finie et dont la théorie élémentaire
est héréditairement indécidable,
alors la théorie élémentaire de $N$ l'est aussi.%
\footnote{Dans \cite{BardakovTolstykh:2007:iaTgF} les auteurs affirment
ce fait avec référence à \cite{Ershov:1980:prkm-rus}.
Nous démontrons ce fait comme le lemme \ref{lemma:19.(6.2)}.}
Ainsi Bardakov et Tolstykh ont démontré que 
la théorie élémentaire de $F$ est héréditairement indécidable,
et donc une partie de la question numéro 4.16 par Mark Sapir dans
\cite{anonymous:2004:Tg40y} est résolue.
Par le même raisonnement, les théories élémentaires de
tous les groupes que nous étudions dans ce travail
sont héréditairement indécidables.

Pour confort du lecteur,
nous présentons dans la section \ref{section:undecidability}
notre version d'une preuve,
basée elle aussi sur un résultat de Mostowski, Robinson, et Tarski
\cite[théorème 9]{MostowskiRobTar:1971:ueua},
que si une structure $S$ de signature finie
interprète l'Arithmétique avec paramètres,
alors la théorie élémentaire de $S$ est héréditairement indécidable.

Les groupes qui font l'objet de notre étude
apparaissent naturellement comme des généralisations
des trois groupes définis par Thompson en 1965 et
habituellement notés $F$, $T$, et $V$.%
\footnote{D'autres lettres ont aussi été utilisées
pour noter ces groupes
(voir \cite{CannonFloydParry:1996:inRTg}).
Il est assez habituel, par exemple,
de noter le groupe $V$ par~$G$.}
Les groupes de Thompson sont exposés en détail dans
\cite{BelkBrown:2005:fdeTgF,CannonFloydParry:1996:inRTg}.
Tous les trois sont infinis et de présentation finie.
Le groupe $V$ était le premier exemple connu d'un groupe simple infini
de présentation finie.
Le groupe $T$ est simple aussi.
Le groupe $F$ se plonge dans $T$, et $T$ se plonge dans $V$.
Le groupe $V$ a été généralisé par Higman \cite{Higman:1974:fpisg}
en une série de groupes $G_{n,r}$, $n=2,3,4,\dotsc$, $r=1,2,3,\dotsc$,
de présentation finie, où $G_{2,1}\cong V$.
Le groupe $G_{n,r}$ de Higman est simple lorsque $n$ est pair ;
lorsque $n$ est impair, le sous-groupe dérivé $[G_{n,r},G_{n,r}]$ est
simple et d'indice $2$ dans $G_{n,r}$.
Kenneth Brown \cite[section 4]{Brown:1987:fpg} a généralisé de la même
façon les groupes $F$ et~$T$.

Les groupes de Thompson ont des représentations
par des permutations affines par morceaux d'un intervalle,
où le groupe $F$ est représenté par
des homéomorphismes par rapport à la topologie habituelle,
et $T$ est représenté par
des homéomorphismes par rapport à la topologie de cercle.
Stein \cite{Stein:1992:gplh} a étudié
trois familles de groupes de telles permutations qui
généralisent respectivement les trois groupes de Thompson.
Afin d'énoncer nos résultats principaux, nous réviserons ici
les définitions de ces familles.

Soient $r$ un nombre réel positif et $\Lambda$ un sous-groupe
du groupe multiplicatif $\mathbb R^*_+$ des nombres réels positifs.
Soit $A$ un sous-groupe additif de $\mathbb R$ contenant $r$ et invariant
sous l'action de $\Lambda$ par multiplication.
Alors définissons $\V(r,\Lambda,A)$ le groupe de toutes les
bijections $x\colon [0;r[\to[0;r[$ qui satisfont
les conditions suivantes :
\begin{enumerate}
\item
	$x$ est affine par morceaux avec un nombre fini des coupures et
	des singularités ;
\item
	$x$ est continue à droite en tout point (au sens habituel) ;
\item
	la pente de chaque partie affine de $x$ est dans $\Lambda$ ;
\item
	tous points de coupure et de singularité de $x$,
	ainsi que leurs images, sont dans~$A$.
\end{enumerate}
La famille des groupes $\V(r,\Lambda,A)$ contiens tous les groupes
de Higman :
pour tout $n=2,3,\dotsc$ et tout $r=1,2,\dotsc$,
$$
G_{n,r}\cong\V(r,\langle n\rangle,\mathbb Z[{\textstyle\frac{1}{n}}]).
$$

Définissons les sous-groupes $\F(r,\Lambda,A)$ et $\T(r,\Lambda,A)$
de $\V(r,\Lambda,A)$ comme suit :
\begin{itemize}
\item
	$\F(r,\Lambda,A)$ est le sous-groupe de tous les éléments
	de $\V(r,\Lambda,A)$ continus par rapport à la topologie habituelle
	de $[0;r[$,
\item
	$\T(r,\Lambda,A)$ est le sous-groupe de tous les éléments de
	$\V(r,\Lambda,A)$ continus par rapport à la topologie
	de cercle sur $[0;r[$
\end{itemize}
(où la topologie de cercle sur $[0;r[$ est la topologie induite
par l'identification naturelle de $[0;r[$ avec
le quotient topologique $[0;r]/\{0,r\}$).
Les groupes $F$, $T$, et $V$ de Thompson sont isomorphes à
$\F(1,\langle 2\rangle,\mathbb Z[1/2])$,
$\T(1,\langle 2\rangle,\mathbb Z[1/2])$,
et $\V(1,\langle 2\rangle,\mathbb Z[1/2])$, respectivement.
Des groupes de la forme $\F(r,\Lambda,A)$ ont été étudiés déjà par
Robert Bieri et Ralph Strebel dans 
\cite{BieriStrebel:pp1985:gPLhrl} (non publié).

Pour le reste nous supposons toujours que $\Lambda\ne\{1\}$.

\begin{TheoremA}
Si\/ $G$ est un sous-groupe de\/ $\V(r,\mathbb R^*_+,\mathbb R)$
tel que
$$
G\cap\F(r,\mathbb R^*_+,\mathbb R)=\F(r,\Lambda,A),
$$
alors\/ $G$ interprète l'Arithmétique\/ $(\mathbb N,+,\times)$
avec paramètres\textup.
\end{TheoremA}


\begin{TheoremB}
Si\/ $\Lambda$ est cyclique\textup,
alors\/ $\F(r,\Lambda,A)$ interprète l'Arithmétique
sans paramètres\textup.
\end{TheoremB}


\begin{TheoremC}
Si\/ $G$ est un groupe comme dans le théorème\/ \textup{A}\textup,
alors la théorie élémentaire de\/ $G$ est
héréditairement indécidable\textup.
\end{TheoremC}

En particulier, tous les groupes de Thompson et de Higman interprètent
l'Arithmétique avec paramètres,
alors que le groupe $F$ de Thompson l'interprète aussi sans paramètres,
et les théories élémentaires de tous ces groupes
sont héréditairement indécidables.

Les théorèmes A et B sont démontrés dans
la section \ref{section:interpretationarithmetic}.
À notre connaissance, l'interprétation construite
dans la preuve du théorème B est entièrement originale.
Le théorème C est démontré dans la section \ref{section:undecidability}
comme un corollaire du théorème A.
Dans l'appendice, nous démontrons que tout élément
du sous-groupe dérivé de $\F(r,\Lambda,A)$ est le produit de deux
commutateurs, et donc que le sous-groupe dérivé est définissable
dans $\F(r,\Lambda,A)$.

L'idée principale de la preuve du théorème A est,
comme dans \cite{BardakovTolstykh:2007:iaTgF},
de trouver dans $G$ un sous-groupe définissable
isomorphe au produit en couronne restreint $\mathbb Z\wr\mathbb Z$,
parce qu'il est connu que ce dernier groupe interprète l'Arithmétique.
Remarquons que, en contraste avec $\mathbb Z\wr\mathbb Z$
et avec les groupes à l'étude,
ni les groupes abéliens, ni les groupes virtuellement abéliens,
ni les groupes libres,
ni les groupes hyperboliques sans torsion
ne peuvent interpréter l'Arithmétique car
leurs théories élémentaires sont tous \emph{stables},
alors que la théorie élémentaire de toute structure
qui interprète l'Arithmétique est «~fortement~» instable.
La \emph{stabilité\/} est une notion fondamentale de la théorie
des modèles, à l'étude de laquelle
\cite{Pillay:1983:ist,Poizat:1985:ctm-fr}
sont d'excellentes introductions.
Les meilleures sources pour l'étude des
\emph{groupes stables}, c'est-à-dire des groupes dont
les théories élémentaires sont stables,
ce sont, à notre avis,
\cite{Poizat:1987:gs-fr,Poizat:2001:sg-eng,Wagner:2000:sg}.
Une démonstration de la stabilité des
groupes abéliens se trouve
dans \cite[théorème 3.1]{Prest:1988:mtm}.
Tout groupe hyperbolique non élémentaire sans torsion est stable
selon un résultat récent de Zlil Sela
\cite{Sela:pp2006:dgg8s}.

\begin{figure}
\includegraphics{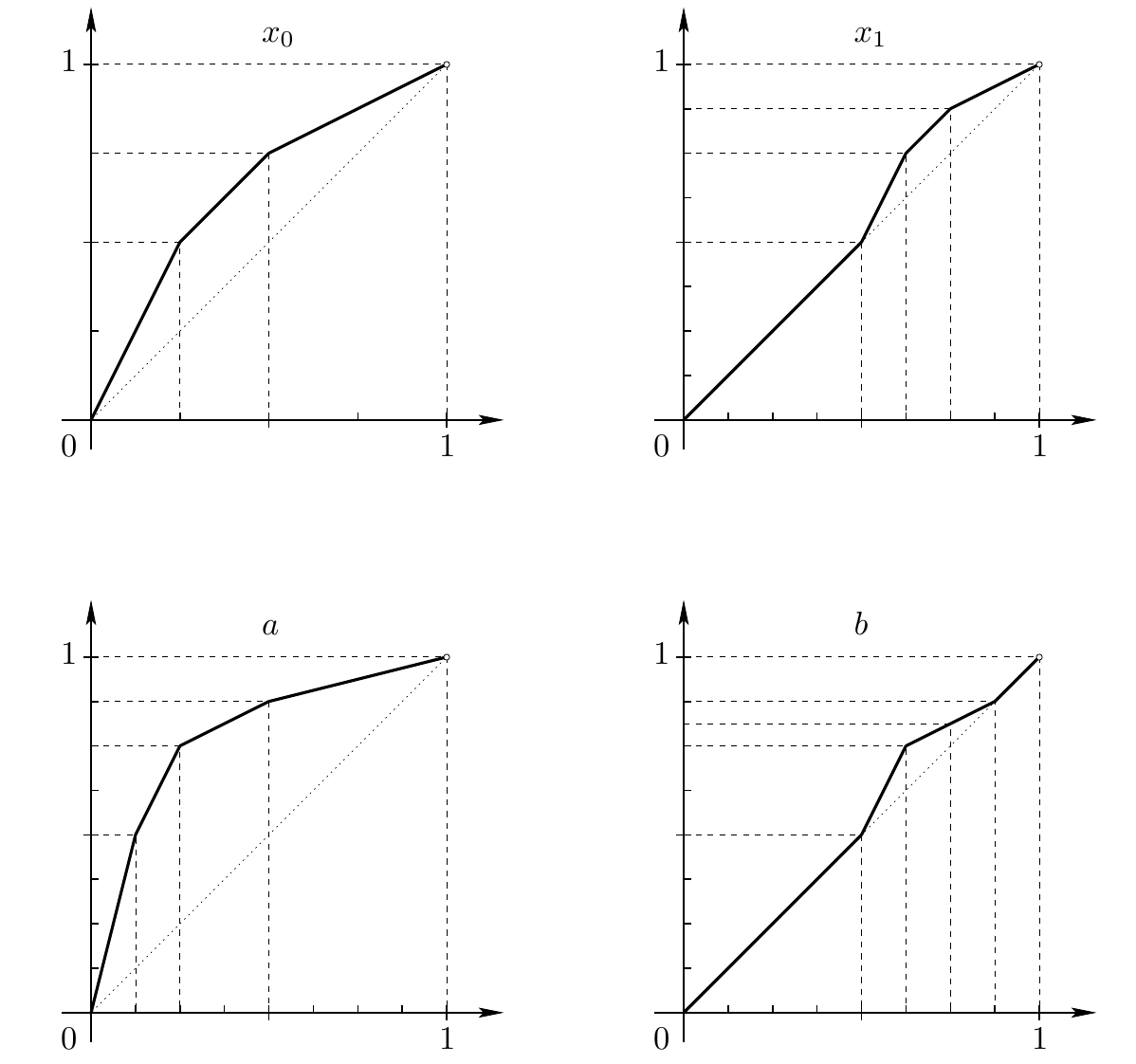}
\caption{Applications $x_0$, $x_1$, $a=x_0^2$,
et $b=x_1x_0^{-1}x_1^{-1}x_0$.}
\label{figure:1}
\end{figure}

Un sous-groupe définissable de $F$ isomorphe à $\mathbb Z\wr\mathbb Z$
a été choisi par Bardakov et Tolstykh
\cite{BardakovTolstykh:2007:iaTgF} comme suit.
Soient $x_0$ et $x_1$ les générateurs «~standards~» de $F$,
et soient $a=x_0^2$, $b=x_1x_0^{-1}x_1^{-1}x_0$
(voir la figure \ref{figure:1}).
On peut vérifier sans grande difficulté que
$\langle a,b\rangle
=\langle b\rangle\wr\langle a\rangle\cong\mathbb Z\wr\mathbb Z$.
Le centralisateur de $x_0$ dans $F$ coïncide avec le sous-groupe
engendré par $x_0$.
En conséquence, le sous-groupe $\langle a\rangle$ est
définissable dans $F$ par une formule avec le paramètre $x_0$.
Puis il est montré que le centralisateur du sous-ensemble
$\{\,a^{-k}ba^{k}\mid k\in\mathbb Z\,\}$ coïncide avec
le sous-groupe $\langle\,a^{-k}ba^{k}\mid k\in\mathbb Z\,\rangle$.
Comme le sous-ensemble $\{\,a^{-k}ba^{k}\mid k\in\mathbb Z\,\}$
est clairement définissable avec les paramètres $x_0$ et $x_1$,
il en est de même pour le sous-groupe
$\langle\,a^{-k}ba^{k}\mid k\in\mathbb Z\,\rangle$.
Donc le sous-groupe $\langle a,b\rangle\cong\mathbb Z\wr\mathbb Z$
est définissable dans $F$ avec paramètres.
Dans la preuve du théorème A, nous suivons une approche similaire
pour le groupe $\F(r,\Lambda,A)$.


\section{Généralités}
\label{section:generalities}

Dans cette section, nous présentons des définitions de base et
quelques faits élémentaires.

\subsection{Permutations et applications affines par morceaux}
\label{subsection:generalities.permutations}

\begin{definition}
Une bijection d'un ensemble sur lui-même est dite une
\emph{permutation\/} de cet ensemble.
Une application $f$ est dite de \emph{permuter\/} un ensemble $S$
si la restriction $f|_S$ est une permutation de~$S$.
\end{definition}


\begin{definition}
Soient $S$ un ensemble et $f$ une bijection de $S$ sur lui-même.
Appelons \emph{support\/} de $f$, noté $\supp(f)$,
le complément dans $S$ de l'ensemble des \emph{points fixes\/} de $f$,
noté $\fix(f)$.
\end{definition}


Par coutume, dans le contexte de l'étude des groupes
de Thompson et de Higman
\emph{toutes les applications agissent à droite}.
Nous adoptons la même convention dans cet article ; par exemple :
$(\alpha)(xy)=((\alpha)x)y$ si
$x$ et $y$ sont des permutations d'un ensemble $S$, et
$\alpha\in S$.

Nous allons écrire $X^f$, ou parfois $(X)f$, pour noter
l'image de l'ensemble $X$ sous l'application~$f$.

Les lemmes \ref{lemma:01.(2.1)}, \ref{lemma:02.(2.2)},
et \ref{lemma:03.(2.3)} sont évidents :

\begin{lemma}
\label{lemma:01.(2.1)}
Deux permutations d'un même ensemble commutent si leurs supports
sont disjoints\textup.
\end{lemma}

\begin{lemma}
\label{lemma:02.(2.2)}
Soient\/ $f$ et\/ $g$ deux permutations d'un même ensemble\textup.
Alors
\begin{align*}
\fix(g^{-1}fg)&=(\fix(f))g,\\
\supp(g^{-1}fg)&=(\supp(f))g.
\end{align*}
\end{lemma}

\begin{lemma}
\label{lemma:03.(2.3)}
Soient\/ $f$ et\/ $g$ deux permutations d'un même ensemble
qui commutent\textup.
Alors\/ $g$ permute chacun des ensembles\/
$\fix(f)$ et\/ $\supp(f)$\textup.
\end{lemma}

\begin{lemma}
\label{lemma:04.(2.4)}
Soient\/ $I$ un intervalle dans\/ $\mathbb R$\textup,
$f$ une permutation croissante de\/ $I$\textup, et\/
$n\in\mathbb Z\setminus\{0\}$\textup.
Alors\/ $\supp(f^n)=\supp(f)$\textup.
\end{lemma}
\begin{proof}
Clairement $\fix(f)\subset\fix(f^n)$ et $\supp(f)\supset\supp(f^n)$.
Considérons $\alpha\in\supp(f)$ arbitrairement choisi.
Sans perte de généralité, supposons que $(\alpha)f>\alpha$.
Alors
$$
\alpha<(\alpha)f<(\alpha)f^2<\dotsb<(\alpha)f^n,
$$
et donc $\alpha\in\supp(f^n)$.
\end{proof}

\begin{lemma}
\label{lemma:05.(2.5)}
Soient\/ $I$ un intervalle compact dans\/ $\mathbb R$\textup,
$f$ une permutation croissante de\/ $I$\textup, et\/
$\alpha\in I$\textup.
Alors
$$
\lim_{n\to+\infty}(\alpha)f^n\in\fix(f).
$$
\end{lemma}
\begin{proof}
Comme $f$ préserve l'ordre, la suite
$((\alpha)f^n)_{n=0,1,\dotsc}$ est monotone, et donc la limite existe
et elle appartient à $I$.
Soit $\beta=\lim_{n\to+\infty}(\alpha)f^n$, alors
$(\beta)f=\beta$ par continuité.
\end{proof}

\begin{definition}
Soient $f$ une application et $\alpha$ un nombre réel.
Disons que $f$ est
\emph{affine à droite de\/ $\alpha$} s'il existe
$\beta>\alpha$ tel que la restriction $f|_{]\alpha;\beta[}$ soit
une application affine $]\alpha;\beta[\to\mathbb R$.
Disons que $f$ est
\emph{affine à gauche de\/ $\alpha$} s'il existe
$\beta<\alpha$ tel que $f|_{]\beta;\alpha[}$ soit
une application affine $]\beta;\alpha[\to\mathbb R$.
\end{definition}

\begin{definition}
Pour tous $\alpha,\beta\in\mathbb R$ tels que $\alpha<\beta$, et
pour toute application $f$ telle que $f|_{]\alpha;\beta[}$ soit une
application affine $]\alpha;\beta[\to\mathbb R$,
nous noterons la pente de $f|_{]\alpha;\beta[}$ par $(\alpha)f^{\prime+}$
et également par $(\beta)f^{\prime-}$.
Disons que $(\alpha)f^{\prime+}$ est
\emph{la pente de\/ $f$ à droite de\/ $\alpha$}, et que
$(\beta)f^{\prime-}$ est
\emph{la pente de\/ $f$ à gauche de\/~$\beta$}.
\end{definition}

On démontre le lemme suivant de la même façon qu'on détermine
la fonction dérivée d'une fonction composée.

\begin{lemma}
\label{lemma:06.(2.6)}
Soient\/ $f$ et\/ $g$ deux applications des sous-ensembles de\/
$\mathbb R$ vers\/ $\mathbb R$\textup, et $\alpha\in\mathbb R$\textup.
\begin{enumerate}
\item
	Si\/ $f$ est affine à droite de\/ $\alpha$ avec\/
	$(\alpha)f^{\prime+}>0$\textup,
	$f$ est continue en\/ $\alpha$ à droite\textup, et\/
	$g$ est affine à droite de\/ $(\alpha)f$\textup, alors
	$$
	(\alpha)(fg)^{\prime+}=(\alpha)f^{\prime+}\cdot((\alpha)f)g^{\prime+}.
	$$
\item
	Si\/ $f$ est affine à gauche de\/ $\alpha$ avec\/
	$(\alpha)f^{\prime-}>0$\textup,
	$f$ est continue en\/ $\alpha$ à gauche\textup, et\/
	$g$ est affine à gauche de\/ $(\alpha)f$\textup, alors
	$$
	(\alpha)(fg)^{\prime-}=(\alpha)f^{\prime-}\cdot((\alpha)f)g^{\prime-}.
	$$
\end{enumerate}
\end{lemma}

\subsection{Groupes à l'étude}
\label{subsection:generalities.groups}

Comme dans l'introduction,
choisissons un nombre réel positif $r$, un sous-groupe
$\Lambda$ du groupe multiplicatif $\mathbb R^*_+$, et
un sous-module $A$ du $\Lambda$-module $\mathbb R$
tel que $r\in A$.
Appelons un tel triple $(r,\Lambda,A)$ \emph{admissible}.
Pour un triple $(r,\Lambda,A)$ admissible,
définissons les groupes
$\F(r,\Lambda,A)$, $\T(r,\Lambda,A)$, et $\V(r,\Lambda,A)$
comme dans l'introduction.
Nous allons traiter les groupes $F$, $T$, et $V$ de Thompson
comme des cas particuliers, donc nous posons
$$
F=\F(1,\langle 2\rangle,\mathbb Z[{\textstyle\frac{1}{2}}]),\qquad
T=\T(1,\langle 2\rangle,\mathbb Z[{\textstyle\frac{1}{2}}]),\qquad
V=\V(1,\langle 2\rangle,\mathbb Z[{\textstyle\frac{1}{2}}]).
$$

\begin{remark}
\label{remark:01.(2.7)}
Tout élément de $\F(r,\mathbb R^*_+,\mathbb R)$ s'étend d'une façon
unique à un homéomorphisme $[0;r]\to[0;r]$
par rapport à la topologie habituelle.
\end{remark}

\begin{remark}
\label{remark:02.(2.8)}
La conjugaison des éléments de $\V(1,\mathbb R^*_+,\mathbb R)$
par l'application linéaire de multiplication par $r$ est
un isomorphisme entre $\V(1,\mathbb R^*_+,\mathbb R)$
et $\V(r,\mathbb R^*_+,\mathbb R)$, qui envoie
$\V(1,\Lambda,Ar^{-1})$ sur $\V(r,\Lambda,A)$,
$\T(1,\Lambda,Ar^{-1})$ sur $\T(r,\Lambda,A)$, et
$\F(1,\Lambda,Ar^{-1})$ sur $\F(r,\Lambda,A)$,
où $Ar^{-1}=\{\,\alpha/r\mid \alpha\in A\,\}$.
\end{remark}

Il est commode de distinguer dans $\F(r,\Lambda,A)$ les sous-semigroupes
suivants :
\begin{itemize}
\item
	notons $\F^\uparrow(r,\Lambda,A)$ le semigroupe de tous les éléments
	$x\in\F(r,\Lambda,A)$ tels que
	$(\alpha)x\ge\alpha$ pour tout $\alpha\in[0;r[$, et 
\item
	notons $\F^\downarrow(r,\Lambda,A)$ le semigroupe de tous les
	$x\in\F(r,\Lambda,A)$ tels que
	$(\alpha)x\le\alpha$ pour tout $\alpha\in[0;r[$.
\end{itemize}

Introduisons d'autres notations utiles :
soit $S\subset[0;r[$, alors
\begin{align*}
\F_S(r,\Lambda,A)
&=\{\,x\in\F(r,\Lambda,A)\mid\supp(x)\subset S\,\},\\
\T_S(r,\Lambda,A)
&=\{\,x\in\T(r,\Lambda,A)\mid\supp(x)\subset S\,\},\\
\V_S(r,\Lambda,A)
&=\{\,x\in\V(r,\Lambda,A)\mid\supp(x)\subset S\,\}.
\end{align*}
On définit de même les semigroupes
$\F^\uparrow_S(r,\Lambda,A)$ et $\F^\downarrow_S(r,\Lambda,A)$.

Pour le reste de cet article, nous fixons $(r,\Lambda,A)$,
et de plus, nous supposons que
\emph{$\Lambda$ n'est pas trivial\/}: $\Lambda\ne\{1\}$.
Afin de simplifier la notation, nous posons :
$$
\F=\F(r,\Lambda,A),\qquad
\T=\T(r,\Lambda,A),\qquad
\V=\V(r,\Lambda,A).
$$
Également, nous allons écrire $\F^\uparrow$ au lieu
de $\F^\uparrow(r,\Lambda,A)$, etc.

\subsection{Théories et modèles}
\label{subsection:generalities.modeltheory}

Dans cet article nous parlons des \emph{structures\/} au sens de
la théorie des modèles
(ou au sens de l'algèbre universelle,
à quelques distinctions linguistiques près).
Lorsque les termes «~formule~», «~énoncé~», et «~théorie~»
sont utilisés dans le sens formel, ils signifient toujours
des formules, des énoncés, et des théories du premier ordre.
Les termes «~théorie~», «~théorie élémentaire~»,
et «~théorie du premier ordre~» seront utilisés comme des synonymes.
«~Un modèle~» et «~une structure~» seront généralement synonymes aussi.
Sauf indication contraire, les formules sont \emph{sans paramètres}.

Une structure $M$ de \emph{signature\/} $\Sigma$,
dite aussi \emph{$\Sigma$-structure\/},
est un \emph{modèle\/}
d'un ensemble d'énoncés $S$, ce qui est noté $M\models S$, si
$M$ \emph{satisfait\/} tout énoncé $\alpha$ de $S$,
ce qui est noté $M\models\alpha$.
Un énoncé $\alpha$ est dit une \emph{conséquence\/}
d'un ensemble d'énoncés $S$ dans une signature $\Sigma$,
ce qui est noté $S\vdash_\Sigma\alpha$
ou $S\vdash\alpha$ au cas où $\Sigma$ est bien comprise, si
tout $\Sigma$-modèle de $S$ est aussi un modèle de $\alpha$.
Un ensemble d'énoncés est \emph{consistant\/}
dans une signature $\Sigma$
s'il a un $\Sigma$-modèle.%
\footnote{La signature $\Sigma$ dans cette définition
est d'importance mineure sauf au cas où $S$ a un modèle vide,
ce qui signifierait en particulier qu'aucun
élément de $S$ ne comporte de symboles de constante.}
Un ensemble d'énoncés est \emph{déductivement clos\/}
dans une signature $\Sigma$
s'il contient toutes ses conséquences dans $\Sigma$.
Une \emph{théorie de signature\/ $\Sigma$}, dite aussi
\emph{$\Sigma$-théorie}, est un ensemble d'énoncés
consistant et déductivement clos dans $\Sigma$.
Si $T$ est une $\Sigma$-théorie, alors $T\vdash_\Sigma\alpha$ équivaut à
$\alpha\in T$.
La \emph{théorie d'une structure\/ $M$}, notée $\Th(M)$,
est l'ensemble de tous les énoncés dans la signature de $M$
satisfaits par $M$.
Une théorie est dite \emph{complète\/}
si elle est la théorie d'une structure.
La classe de tous les $\Sigma$-modèles d'un ensemble
de $\Sigma$-énoncés $S$ est notée $\Mod_\Sigma(S)$.

Nous allons utiliser implicitement \emph{le théorème de compacité},
qui assure qu'une conséquence d'un ensemble d'énoncés est toujours une
conséquence, dans la même signature, d'une de ses parties finies.
Nous recommandons un des ouvrages
\cite{Hodges:1993:mt,Hodges:1997:smt,
Poizat:1985:ctm-fr,Rothmaler:2000:imt-eng}
pour références sur des résultats généraux
de la théorie de modèles.

Si $M$ est une structure,
un \emph{ensemble définissable dans\/} $M$ est en général
une partie de $M^n$, où $n\in\mathbb N$,
définissable par une formule du premier ordre dans la
signature de $M$, et éventuellement avec
des \emph{paramètres\/} extraits de $M$.
Les \emph{paramètres\/} sont de nouveaux symboles de constante rajoutés
dans le langage et interprétés par des éléments de $M$.
(D'habitude, pour nommer un élément $a\in M$,
on utilise $a$ lui-même comme un paramètre.)
Par exemple, dans un groupe, le centralisateur de tout élément $g$ est
définissable par la formule $\phi(x)=\ulcorner gx=xg\urcorner$
avec $g$ comme un paramètre,
mais en général il n'y a aucune raison pour que ce centralisateur
soit définissable par une formule sans paramètres
dans la signature de groupe pur,
qui ne comporte qu'un seul symbole de fonction binaire
$\ulcorner\cdot\urcorner$ pour noter l'opération
de groupe $(x,y)\mapsto x\cdot y$.
Pour préciser, un ensemble est dit
«~définissable \emph{avec\/}~» ou «~\emph{sans\/} paramètres~»
selon si des paramètres sont permis ou pas dans sa définition.
Par contre, nous n'allons pas préciser la structure dans laquelle un
ensemble donné sera définissable, sauf s'il y a plusieurs
choix également naturels dans le contexte.
Bien entendu, on peut parler aussi de la définissabilité des relations
et des opérations.

Si $f$ est une application $A\to B$, et que $n$ est un entier positif,
nous notons $f^{\underline{n}}$ l'application $A^n\to B^n$
induite par $f$.
Nous allons légèrement abuser la notation en supposant que si l'ensemble
de départ de $f$ est une partie de $A^m$, alors l'ensemble de départ
de $f^{\underline{n}}$ est identifié naturellement avec
une partie de $A^{mn}$.
Nous appelons la \emph{$f$-préimage\/} d'un ensemble donné
sa préimage sous $f^{\underline{n}}$
au cas où le choix de $n$ est évident
(donc pas nécessairement sous $f$ elle-même).

Considérons deux structures,
$M$ de signature $\Sigma$ et $N$ de signature $\Gamma$.

\begin{definition}
Nous appelons une \emph{interprétation avec paramètres de\/ $M$
dans\/ $N$} une paire $(n,f)$ où
$n\in\mathbb N$ et où $f$ est une application surjective
d'une partie de $N^n$ sur $M$,
telle que pour tout ensemble $X$ définissable dans $M$ sans paramètres,
la $f$-préimage de $X$ soit définissable (dans $N$)
avec (éventuellement) paramètres.
Une interprétation $(n,f)$ avec paramètres est dite une
\emph{interprétation sans paramètres\/} si
la $f$-préimage de tout ensemble définissable
sans paramètres est, elle aussi, définissable sans paramètres.
\end{definition}

Voir \cite[chapitre 5]{Hodges:1993:mt} pour des explications
détaillées de l'interprétabilité et des notions liées.

Dans ce qui suit, conformément à l'usage en théorie des modèles,
les termes «~définissable~», «~$0$-définissable~», 
«~interprétation~», et «~$0$-in\-ter\-pré\-ta\-tion~»
seront utilisés pour noter
«~définissable avec paramètres~»,
«~définissable sans paramètres~»,
«~interprétation avec paramètres~», et
«~interprétation sans paramètres~»,
respectivement.
Aussi, comme dans notre cas la valeur de $n$ pour une interprétation
$(n,f)$ à l'étude sera souvent soit bien comprise, soit peu importante,
nous allons simplifier la notation et appeler
$f$ elle-même une interprétation.
Pour noter que $(n,f)$ est une interprétation de $M$ dans $N$,
nous allons écrire
soit $(n,f)\colon M\rightsquigarrow N$,
soit $f\colon M\rightsquigarrow N$.%
\footnote{La flèche ici se dirige dans le sens inverse
de la notation de \cite{AhlbrandtZiegler:1986:qfatct}.}
Pour noter que $f$ est une $0$-interprétation de $M$ dans $N$,
nous écrivons
$f\colon M\stackrel{\scriptscriptstyle\varnothing}{\rightsquigarrow}N$.
Pour noter que $M$ est interprétable ou $0$-interprétable
dans $N$, nous écrivons $M\rightsquigarrow N$ ou
$M\stackrel{\scriptscriptstyle\varnothing}{\rightsquigarrow}N$,
respectivement.

\begin{remark}
\label{remark:03.(2.9)}
Si $f\colon M\rightsquigarrow N$, alors
la $f$-préimage de tout ensemble définissable dans $M$ est,
elle aussi, définissable (dans $N$).
\end{remark}

\begin{remark}
\label{remark:04.(2.10)}
Soient $n\in\mathbb N$ et $B\subset N^n$.
Alors une application surjective $f$ de $B$ sur $M$ est une interprétation
de $M$ dans $N$ si et seulement si
\begin{enumerate}
\item
	l'ensemble de départ $B$ est définissable,
\item
	la relation d'équivalence sur $B$ induite par $f$
	(le \emph{noyau\/} de $f$) est définissable, et
\item
	pour toute relation, opération, et constante de la structure $M$
	(nommée par un symbole de $\Sigma$),
	la $f$-préimage de son graphe est définissable.
\end{enumerate}
L'application $f$ est une $0$-interprétation si et seulement si
tous ces ensembles sont $0$-définissables.
\end{remark}

\begin{remark}
\label{remark:05.(2.11)}
Si $L$, $M$, et $N$ sont trois structures, et que
$(m,f)\colon L\rightsquigarrow M$ et $(n,g)\colon M\rightsquigarrow N$,
alors $(mn,g^{\underline{m}}f)\colon L\rightsquigarrow N$.
Si de plus $f$ et $g$ sont $0$-interprétations, alors
$g^{\underline{m}}f$ l'est aussi.
\end{remark}

\begin{definition}
[voir \cite{AhlbrandtZiegler:1986:qfatct} et
{\cite[section 5.4(c)]{Hodges:1993:mt}}]
Deux structures $M$ et $N$ sont dites \emph{bi-interprétables\/}
s'il existe deux interprétations
$(m,f)\colon M\rightsquigarrow N$ et $(n,g)\colon N\rightsquigarrow M$
telles que
l'application $g^{\underline{m}}f$ soit définissable dans $M$,
et que $f^{\underline{n}}g$ soit définissable dans $N$.
Les interprétations $(m,f)$ et $(n,g)$ dans ce cas sont dites
\emph{bi-interprétations}.
\end{definition}

\subsection{Décidabilité}
\label{subsection:generalities.decidability}

Soit $A$ un ensemble fini vu comme un alphabet,
et notons $A^*$ l'ensemble de tous les mots finis dans $A$.
Nous appelons un ensemble $X\subset A^*$ \emph{récursif\/}
ou bien \emph{décidable\/}
s'il existe un \emph{algorithme\/} qui pour toute donnée $w\in A^*$
répond à la question si $w\in X$.
Une application $f\colon X\to A^*$ est dite \emph{calculable\/}
s'il existe un algorithme qui calcule $f(w)$ pour toute donnée $w\in X$,
et qui ne s'arrête jamais pour toute donnée $w\notin X$.

Normalement on dit qu'un ensemble est récursif ou non récursif,
tandis que une théorie est décidable ou indécidable.

Dans le reste de cette section, soit $\Sigma$
une signature finie arbitraire.

\begin{definition}
Une $\Sigma$-théorie $T$
est dite \emph{essentiellement indécidable\/} si
toute $\Sigma$-théorie contenant $T$ est indécidable.
\end{definition}

\begin{remark}
\label{remark:06.(2.12)}
Si $T$ est une théorie de signature finie,
et que une partie de $T$ est une théorie essentiellement indécidable
(de signature éventuellement plus petite),
alors $T$ est indécidable, et même essentiellement indécidable.
\end{remark}

\begin{definition}
Une $\Sigma$-théorie $T$
est dite \emph{héréditairement indécidable\/} si
toute $\Sigma$-sous-théorie de $T$ est indécidable.
\end{definition}

\begin{lemma}[{\cite[théorème 6]{Tarski:1971:gmpu}}]
\label{lemma:07.(2.13)}
Si\/ $T$ est une théorie de signature finie\textup,
et que\/ $T$ a une sous-théorie\/
\textup(de signature éventuellement plus petite\/\textup)
essentiellement indécidable et finiment axiomatisée\textup,
alors\/ $T$ est héréditairement indécidable\textup.
\end{lemma}
\begin{proof}
Notons $\Sigma$ la signature de $T$.
Soit $S$ une sous-théorie de $T$ qui est essentiellement indécidable et
finiment axiomatisée.
Choisissons un énoncé $\theta\in S$ qui axiomatise~$S$.

Soit par l'absurde $U$ une $\Sigma$-sous-théorie de $T$
qui est décidable.
Soit $R$ la $\Sigma$-théorie axiomatisée (engendrée) par $U\cup S$.
Alors
$$
R=\{\,\alpha\mid U,\theta\vdash_\Sigma\alpha\,\}
=\{\,\alpha\mid\ulcorner\theta\rightarrow\alpha\urcorner\in U\,\},
$$
et donc $R$ est décidable puisque $U$ l'est.
Cela contredit l'indécidabilité essentielle de $S$
(voir la remarque~\ref{remark:06.(2.12)}).
\end{proof}


\section{Deux lemmes}
\label{section:twolemmas}

Dans cette section, nous démontrerons deux lemmes techniques à propos
des homéomorphismes affines par morceaux d'un intervalle.
Ces lemmes seront essentiels pour la preuve du théorème A, à savoir,
pour démontrer que certains centralisateurs sont préservés en passant
de $\F$ à $\V$, et ainsi pour pouvoir passer d'une interprétation
de l'Arithmétique dans $\F$ à ses interprétations dans $\V$ et dans~$\T$.

\begin{lemma}
\label{lemma:08.(3.1)}
Soit\/ $r\in\mathbb R^*_+$\textup.
Soit\/ $z$ un homéomorphisme\/ $[0;r[\to[0;r[$ tel que\/
$(\alpha)z>\alpha$ pour tout\/ $\alpha\in]0;r[$\textup.
Soit\/ $f$ une permutation de\/ $[0;r[$
telle que\/ \textup:
\begin{enumerate}
\item
	$f$ commute avec\/ $z$\textup,
\item
	$f$ soit continue\/ \textup(à droite\/\textup) en\/ $0$\textup,
\item
	$f$ ait un nombre fini des points de discontinuité\textup.
\end{enumerate}
Alors\/ $f$ est continue\/
\textup(par rapport à la topologie habituelle\/\textup)\textup.
\end{lemma}
\begin{proof}
Supposons que $f$ ne soit pas continue.
Alors soit $\alpha$ le plus petit élément de $[0;r[$
où $f$ n'est pas continue.
Comme $f$ est continue en $0$, $\alpha\in]0;r[$.
Donc $(\alpha)z^{-1}<\alpha$ et $f$ est continue en $(\alpha)z^{-1}$.
On conclut que $f$ est continue en $\alpha$ parce que
$f=z^{-1}fz$, où $z^{-1}$ et $z$ sont continues partout,
et $f$ est continue en $(\alpha)z^{-1}$.
Cela donne une contradiction.
\end{proof}

\begin{lemma}
\label{lemma:09.(3.2)}
Soit\/ $r\in\mathbb R^*_+$\textup.
Soit\/ $Z$ un ensemble d'homéomorphismes\/ $[0;r[\to[0;r[$
tel que\/ \textup:
\begin{enumerate}
\item
	$(\alpha)z\ge\alpha$\/ pour tout\/ $z\in Z$ et tout\/
	$\alpha\in[0;r[$\textup,
\item
	$\supp(z)$ soit un intervalle pour tout\/ $z\in Z$\textup,
\item
	$\bigcup_{z\in Z}\supp(z)$ soit dense dans\/ $]0;r[$\textup.
\end{enumerate}
Soit\/ $f$ une permutation de\/ $[0;r[$ telle que\/ \textup:
\begin{enumerate}
\item
	$f$ commute avec tout\/ $z\in Z$\textup,
\item
	$f$ soit continue à droite en tout point de\/ $[0;r[$\textup,
\item
	$f$ ait un nombre fini des points de discontinuité\textup.
\end{enumerate}
Alors\/ $f$ est continue\textup.
\end{lemma}
\begin{proof}
Démontrons d'abord que $f$ est croissante
sur chacun de ses intervalles de continuité.

Disons qu'un intervalle $I$ est \emph{fermé à droite\/}
si $\sup I\in I$,
et qu'il est \emph{fermé à gauche\/}
si $\inf I\in I$.
Si un intervalle n'est pas fermé à droite ou gauche, disons qu'il
\emph{y\/} est \emph{ouvert}.

Comme le nombre de points de discontinuité de $f$ est fini,
le nombre de ses intervalles de continuité maximaux est fini aussi.
Comme $f$ est continue à droite partout
dans son ensemble de départ $[0;r[$,
tout intervalle de continuité maximal est fermé à gauche et
ouvert à droite.
Alors l'image sous $f$ de tout intervalle $I$ de continuité maximal
est fermé à gauche si $f$ est croissante sur $I$,
et fermé à droite si $f$ est décroissante sur~$I$.

Évidemment $[0;r[$ est la réunion disjointe des images
des intervalles de continuité maximaux de $f$.
Comme toute telle image est un intervalle soit ouvert à droite
soit ouvert à gauche, la seule possibilité est
qu'ils sont tous fermés à gauche et ouverts à droite.
Donc $f$ est croissante sur chacun de ses intervalles de continuité.

Supposons maintenant que $f$ ne soit pas continue.

D'après le lemme \ref{lemma:03.(2.3)}, pour tout $z\in Z$,
$f$ permute l'ensemble $\supp(z)$.
Démontrons que $f$ est continue sur chacun de ces intervalles.
Soit par l'absurde $z\in Z$ tel que
$f$ ne soit pas continue sur $\supp(z)$.
Soit $\gamma$ le plus petit élément de $\supp(z)$ où
$f$ n'est pas continue.
Alors $f$ est continue en $\gamma$ parce que
$f=z^{-1}fz$, où $z$ est un homéomorphisme,
et $f$ est continue en $(\gamma)z^{-1}$, puisque $(\gamma)z^{-1}<\gamma$.
Cela donne une contradiction.

Pour tout $z\in Z$, $f$ est croissante sur $\supp(z)$
puisqu'elle y est continue.
Alors pour tout $z\in Z$,
$$
\lim_{\alpha\to(\sup\supp(z))^{-}}(\alpha)f=\sup\supp(z).
$$
Par continuité à droite, $(\inf\supp(z))f=\inf\supp(z)$
pour tout $z\in Z$.

Soient $S=\bigcup_{z\in Z}\supp(z)$ et
$L=\{\,\inf\supp(z)\mid z\in Z\,\}$.
Alors $S$ est ouvert et dense dans $[0;r[$, $f|_S$ est continue,
et $f|_L=\id_L$.

Considérons $\alpha\in[0;r[\setminus S$ arbitraire.
Il est facile de voir que pour tout $\beta\in]\alpha;r[$,
$L\cap[\alpha;\beta[$ n'est pas vide.
Donc $(\alpha)f=\alpha$ par continuité à droite en $\alpha$.
Nous avons démontré que
$f|_{[0;r[\setminus S}=\id_{[0;r[\setminus S}$.

L'application $f$ est croissante.
En effet, si $f$ est une application d'un ensemble linéairement ordonné
vers lui-même, et cet ensemble est recouvert par des intervalles tels
que $f$ envoie chacun d'entre eux vers lui-même d'une façon
strictement croissante,
alors $f$ est strictement croissante.
Dans notre cas, nous avons :
$$
[0;r[=\bigcup_{z\in Z}\supp(z)
\cup\bigcup_{\alpha\in[0;r[\setminus S}[\alpha;\alpha].
$$

Alors $f$ est continue parce qu'elle est une surjection croissante d'un
sous-ensemble de $\mathbb R$ sur un intervalle de~$\mathbb R$.
\end{proof}

Ces deux derniers lemmes déjà suffisent pour démontrer,
à partir du résultat de Bardakov et Tolstykh,
que le sous-groupe définissable de $F$ isomorphe à $\mathbb Z\wr\mathbb Z$
utilisé dans \cite{BardakovTolstykh:2007:iaTgF} est définissable
dans $T$ et $V$ aussi.

\begin{proposition}
\label{proposition:01.(3.3)}
Soient\/ $a$ et\/ $b$ les éléments de\/ $F$ montrés sur
la figure\/ \textup{\ref{figure:1}.}
Alors
$$
\langle a,b\rangle
=\langle b\rangle\wr\langle a\rangle\cong\mathbb Z\wr\mathbb Z,
$$
et le sous-groupe\/ $\langle a,b\rangle$ est définissable
dans\/ $F$\textup, dans\/ $T$\textup, et dans\/ $V$
par la même formule du premier ordre avec paramètres\textup.
\end{proposition}
\begin{proof}
Nous considérons les mêmes éléments
$x_0$, $x_1$, $a=x_0^2$, et $b=x_1x_0^{-1}x_1^{-1}x_0$ de $F$
que dans \cite{BardakovTolstykh:2007:iaTgF},
voir la figure \ref{figure:1}.
Il est démontré dans \cite{BardakovTolstykh:2007:iaTgF} que :
\begin{enumerate}
\item
	$\langle a,b\rangle
	=\langle b\rangle\wr\langle a\rangle\cong\mathbb Z\wr\mathbb Z$,
\item
	$C_F(x_0)=\langle x_0\rangle$,
\item
	$C_F(\{\,a^{-k}ba^{k}\mid k\in\mathbb Z\,\})
	=\langle\,a^{-k}ba^{k}\mid k\in\mathbb Z\,\rangle$.
\end{enumerate}
En particulier, $\langle a,b\rangle$ est le produit semi-direct de
$C_F(\{\,a^{-k}ba^{k}\mid k\in\mathbb Z\,\})$ et de $\langle a\rangle$,
où $\langle a\rangle=\{\,x^2\mid x\in C_F(x_0)\,\}$.

Notons $\alpha_k=2^{-1+2k}$ pour $k=-1,-2,-3,\dotsc$, et
$\alpha_k=1-2^{-1-2k}$ pour $k=0,1,2,\dotsc$.
Alors
$$
0<\dotsb<\alpha_{-2}<\alpha_{-1}<\alpha_0<\alpha_1<\alpha_2<\dotsb<1.
$$

Un calcul direct facilité par le lemme \ref{lemma:02.(2.2)} montre que :
\begin{enumerate}
\item
	$\supp(x_0)=]0;1[$ ;
\item
	$(\alpha)x_0\ge\alpha$ pour tout $\alpha\in[0;1[$ ;
\item
	$\supp(a^{-k}ba^k)=]\alpha_k;\alpha_{k+1}[$ pour tout $k\in\mathbb Z$ ;
\item
	$(\alpha)a^{-k}ba^k\ge\alpha$
	pour tout $\alpha\in[0;1[$ et tout $k\in\mathbb Z$.
\end{enumerate}

Si on prend $x_0$ comme $z$ dans le lemme \ref{lemma:08.(3.1)},
ce lemme permet de conclure que tout élément de $C_V(x_0)$ est continu.
Si on pose $Z=\{\,a^{-k}ba^{k}\mid k\in\mathbb Z\,\}$
dans le lemme \ref{lemma:09.(3.2)},
ce lemme montre que tout élément de
$C_V(\{\,a^{-k}ba^{k}\mid k\in\mathbb Z\,\})$ est continu.
Comme un élément de $V$ appartient à $F$ si et seulement si il
est continu, on conclut que
le centralisateur de l'élément $x_0$ et le centralisateur de
l'ensemble $\{\,a^{-k}ba^{k}\mid k\in\mathbb Z\,\}$
sont conservés en passant de $F$ à $V$.
Donc $\langle a,b\rangle$ est définissable dans $F$, $T$, et $V$
par la même formule du premier ordre avec paramètres.
\end{proof}


\section{Copies définissables de $\mathbb Z\wr\mathbb Z$}
\label{section:copiesZwrZ}

Dans cette section, nous démontrons que tous les groupes
$\F$, $\T$, et $\V$ ont de sous-groupes définissables isomorphes
à $\mathbb Z\wr\mathbb Z$.

\begin{lemma}
\label{lemma:10.(4.1)}
Soient\/ $\alpha,\beta\in[0;r]\cap A$ tels que\/ $\alpha<\beta$\textup,
et soit\/ $I=]\alpha;\beta[$\textup.
Soit\/ $x\in\F_I$ tel que\/
$\fix(x)\cap I\cap A=\varnothing$\textup.
Soit\/ $\phi\colon \F_I\to\Lambda$
l'application\/ $y\mapsto(\alpha)y^{\prime+}$\textup.
Soit\/ $C$ le centralisateur de\/ $x$ dans\/
$\F_I$\textup.
Alors\/ $\phi$ est un homomorphisme\textup, et
sa restriction à\/ $C$ est injective\textup.
De même pour\/
$\phi\colon \F_I\to\Lambda,\ y\mapsto(\beta)y^{\prime-}$\textup.
\end{lemma}
\begin{proof}
Considérons seulement le cas de 
$\phi\colon \F_I\to\Lambda,\ y\mapsto(\alpha)y^{\prime+}$, parce que
le cas de $y\mapsto(\beta)y^{\prime-}$ est analogue.
Il est facile de vérifier que $\phi$ est un homomorphisme
(voir le lemme \ref{lemma:06.(2.6)}).
Supposons que $\phi$ ne soit pas injectif sur~$C$.

Soit $y\in C$ tel que $\phi(y)=1$ mais $y\ne\id$.
Soit $\gamma\in]\alpha;\beta[$ tel que
$$
y|_{[0;\gamma]}=\id_{[0;\gamma]}\quad\text{mais}\quad
(\gamma)y^{\prime+}\ne1.
$$
Alors $\gamma\in A$, et donc $(\gamma)x\ne\gamma$.
Sans perte de généralité, supposons que $(\gamma)x>\gamma$,
car sinon, alors $(\gamma)x^{-1}>\gamma$, et on peut utiliser
$x^{-1}$ au lieu de $x$.
Alors
$$
[0;(\gamma)x]=[0;\gamma]^x\subset\fix(y),
$$
voir le lemme \ref{lemma:03.(2.3)},
et donc $(\gamma)y^{\prime+}=1$.
Cela donne une contradiction.
\end{proof}


\begin{lemma}
\label{lemma:11.(4.2)}
Le centralisateur\/ $C$ dans le lemme\/ \textup{\ref{lemma:10.(4.1)}}
est cyclique\textup.
\end{lemma}

Ce lemme résulte de la description des centralisateurs
dans $\F(r,\mathbb R^*_+,\mathbb R)$
obtenu par Matthew Brin et Craig Squier
\cite{BrinSquier:2001:pcrcgplhrl},
mais pour confort du lecteur nous préférons de fournir
une preuve autonome.

\begin{proof}[Démonstration du lemme\/ \textup{\ref{lemma:11.(4.2)}}]
Notons tout d'abord que si $\Lambda$ est cyclique lui-même,
la conclusion de ce lemme est un corollaire évident du
lemme~\ref{lemma:10.(4.1)}.

Soient $\alpha$, $\beta$, $I$, $x$, $\phi$, et $C$
tels comme dans le lemme \ref{lemma:10.(4.1)}.
Sans perte de généralité, supposons que $(\alpha)x^{\prime+}>1$.

Nous allons utiliser le fait que un sous-groupe multiplicatif de
$\mathbb R^*_+$ est soit cyclique (le sous-groupe trivial y compris),
soit dense dans $\mathbb R^*_+$ par rapport à la topologie habituelle.

Soit $\Gamma$ l'image du groupe $C$ sous l'homomorphisme
$\phi\colon \F_I\to\Lambda$.
D'après le lemme \ref{lemma:10.(4.1)}, $\phi$ est injectif,
et donc $C\cong\Gamma$.
Il nous reste à montrer que $\Gamma$ n'est pas dense
dans $\mathbb R^*_+$.

Notons que $\fix(x)\cap I$ est
un ensemble fini, et que
$$
\fix(y)\cap I=\fix(x)\cap I\quad\text{pour tout}\quad
y\in C\setminus\{\id\}.
$$
En effet, il est clair que $\fix(x)\cap I$ est fini parce que
$\fix(x)\cap I\cap A$ est vide et $A$ est dense dans $\mathbb R$.
Considérons maintenant un $y\in C\setminus\{\id\}$.
Par l'injectivité de $\phi$, $(\alpha)y^{\prime+}\ne1$.
Si $\fix(y)\cap I\cap A$ était non vide, alors il aurait
le plus petit élément $\gamma$,
et ce $\gamma$ serait fixé par $x$ : $(\gamma)x=x$.
Cela serait une contradiction avec
$\fix(x)\cap I\cap A=\varnothing$,
donc $\fix(y)\cap I\cap A=\varnothing$
et $\fix(y)\cap I$ est fini.
Comme $x$ et $y$ commutent et que chacune permute $I$,
on voit que $x$ permute $\fix(y)\cap I$,
et que $y$ permute $\fix(x)\cap I$.
Comme ces ensembles sont finis, et que
$x$ et $y$ préservent l'ordre, on conclut que
$\fix(y)\cap I\subset\fix(x)$,
$\fix(x)\cap I\subset\fix(y)$,
et donc $\fix(y)\cap I=\fix(x)\cap I$.

Notons
$$
\beta_0=\min((\fix(x)\cap I)\cup\{\beta\}) ;
$$
il existe mais il n'appartient pas à $A$ sauf si $\beta_0=\beta$.
Alors $(\gamma)x>\gamma$ pour tout $\gamma\in]\alpha;\beta_0[$,
puisque $(\alpha)x^{\prime+}>1$.

Choisissons $\alpha_1,\beta_1\in]\alpha;\beta_0[$ tels que
$x^{-1}$ soit affine sur
$[\alpha;\alpha_1]$ et $x$ soit affine sur $[\beta_1;\beta_0[$.
Alors $x$ est aussi affine sur
$[\alpha;(\alpha_1)x^{-1}]$,
et $x^{-1}$ est affine sur $[(\beta_1)x;\beta_0[$.

Démontrons maintenant que si
$y\in C$ et $(\alpha)y^{\prime+}>1$, alors $y$ est affine sur
$[\alpha;(\alpha_1)y^{-1}]$ et sur $[\beta_1;\beta_0[$.
Considérons un tel $y$.
Donc $(\gamma)y>\gamma$ pour tout $\gamma\in]\alpha;\beta_0[$.
Alors pour tout $\gamma\in]\alpha;\alpha_1]$,
$$
x^{-1}|_{[\alpha;\alpha_1]}\cdot y^{-1}|_{[\alpha;(\gamma)x^{-1}]}\cdot
x|_{[\alpha;(\alpha_1)x^{-1}]}
=y^{-1}|_{[\alpha;\gamma]},
$$
et pour tout $\gamma\in[\beta_1;\beta_0[$,
$$
x|_{[\beta_1;\beta_0[}\cdot y|_{[(\gamma)x;\beta_0[}\cdot
x^{-1}|_{[(\beta_1)x;\beta_0[}
=y|_{[\gamma;\beta_0[}.
$$
Ces égalités évidentes impliquent que :
\begin{enumerate}
\item
	pour tout $\gamma\in]\alpha;\alpha_1]$,
	$y^{-1}$ est affine sur $[\alpha;\gamma]$ dès qu'elle
	est affine sur $[\alpha;(\gamma)x^{-1}]$,
\item
	pour tout $\gamma\in[\beta_1;\beta_0[$,
	$y$ est affine sur $[\gamma;\beta_0[$ dès qu'elle
	l'est sur $[(\gamma)x;\beta_0[$.
\end{enumerate}
Ceci n'est possible que si $y^{-1}$ est affine sur $[\alpha;\alpha_1]$,
et $y$ est affine sur $[\beta_1;\beta_0[$.

Clairement $\lim_{n\to+\infty}(\gamma)x^{n}=\beta_0$
pour tout $\gamma\in]\alpha;\beta_0[$
(voir la preuve du lemme \ref{lemma:05.(2.5)}).
Choisissons un entier $n$ tel que
$$
(\alpha_1)x^{n}>\beta_1,\quad\text{et donc}\quad(\beta_1)x^{-n}<\alpha_1.
$$
Notons $p=(\beta_1)(x^{-n})^{\prime+}$.
Alors $p^{-1}=((\beta_1)x^{-n})(x^{n})^{\prime+}$.
Choisissons $\gamma$ dans $((\beta_1)x^{-n};\alpha_1]$
tel que $x^{n}$ soit affine sur
$[(\beta_1)x^{-n};\gamma]$ (avec la pente $p^{-1}$).

Supposons que $y$ soit un élément de $C$ tel que
$$
1<(y)\phi<\frac{\gamma-\alpha}{(\beta_1)x^{-n}-\alpha}.
$$
Alors $y^{-1}$ est affine sur
$[\alpha;\alpha_1]$, $y$ est affine sur $[\beta_1;\beta_0[$,
et $(\beta_0)y^{\prime-}<1$
(parce que $(\gamma)y>\gamma$ pour tout $\gamma\in]\alpha;\beta_0[$).
Notons $q=(\alpha)y^{\prime+}=(y)\phi$.
Alors
$$
(\beta_1)x^{-n}<\alpha+(\gamma-\alpha)q^{-1}=(\gamma)y^{-1}
\le\alpha+(\alpha_1-\alpha)q^{-1}=(\alpha_1)y^{-1},
$$
et par conséquent $(\beta_1)x^{-n}<(\beta_1)x^{-n}y<\gamma$.
Comme $y=x^{-n}yx^{n}$, on calcule
\begin{align*}
(\beta_1)y^{\prime+}
&=(\beta_1)(x^{-n})^{\prime+}
\cdot((\beta_1)x^{-n})y^{\prime+}
\cdot((\beta_1)x^{-n}y)(x^{n})^{\prime+}\\
&=pqp^{-1}=q>1.
\end{align*}
Comme $(\beta_1)y^{\prime+}=(\beta_0)y^{\prime-}<1$,
cela donne une contradiction, ce qui signifie que
$$
\Gamma\cap\Bigl]1;\frac{\gamma-\alpha}{(\beta_1)x^{-n}-\alpha}\Bigr[
=\varnothing.
$$
\end{proof}

\begin{lemma}
\label{lemma:12.(4.3)}
Soient\/ $\alpha,\beta\in[0;r]\cap A$ tels que\/ $\alpha<\beta$\textup.
Soient\/ $p,q\in\Lambda$ tels que\/ $p>1>q$\textup.
Alors il existe\/ $x\in\F^\uparrow$ tel que\/
$\supp(x)=]\alpha;\beta[$\textup,
$(\alpha)x^{\prime+}=p$\textup, et que\/ $(\beta)x^{\prime-}=q$\textup.
\end{lemma}
\begin{proof}
Soit $s$ un élément de $\Lambda$ tel que $(2+p+q)s\le1$.
Notons $l=\beta-\alpha$.
Faisons deux partitions de l'intervalle $[\alpha;\beta]$ :
la première --- en sous-intervalles de longueurs
$sl$, $qsl$, $(1-(2+p+q)s)l$, $psl$, et $sl$,
dans cet ordre,
et la deuxième --- en sous-intervalles de longueurs
$psl$, $sl$, $(1-(2+p+q)s)l$, $sl$, et $qsl$, dans cet ordre.
Soit $x$ l'application continue $[0;r]\to[0;r]$ qui est
l'identité sur $[0;\alpha]\cup[\beta;r]$, et 
qui envoie chaque intervalle de la première partition de $[\alpha;\beta]$
de la façon affine sur l'intervalle correspondante de la deuxième.
Il est facile de vérifier que
$\supp(x)=]\alpha;\beta[$,
$(\alpha)x^{\prime+}=p$, $(\beta)x^{\prime-}=q$,
et que $x|_{[0;r[}\in\F^\uparrow$.
\end{proof}

\begin{figure}
\includegraphics{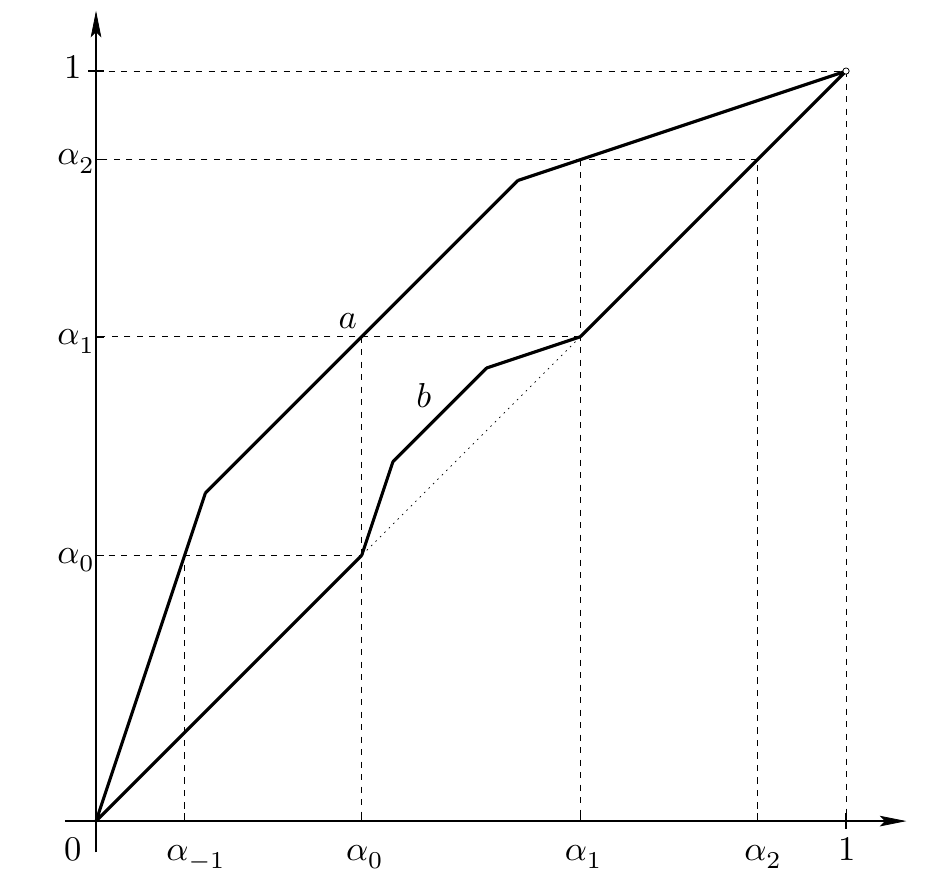}
\caption{Un exemple des applications $a$ et $b$.}
\label{figure:2}
\end{figure}

Choisissons $a\in\F^\uparrow$ tel que $\supp(a)=]0;r[$.
Choisissons $\alpha_0\in]0;r[\cap A$ arbitrairement.
Pour tout $k\in\mathbb Z$, définissons $\alpha_k=(\alpha_0)a^k$.
Notons que
\begin{gather*}
0<\dotsb<\alpha_{-2}<\alpha_{-1}<\alpha_0<\alpha_1<\alpha_2<\dotsb<r,\\
\text{et que}\qquad
\lim_{n\to-\infty}\alpha_n=0,\quad\lim_{n\to+\infty}\alpha_n=r
\end{gather*}
(voir le lemme \ref{lemma:05.(2.5)} et
la remarque \ref{remark:01.(2.7)}).
Choisissons $b\in\F^\uparrow$ tel que
$\supp(b)=]\alpha_0;\alpha_1[$
(voir la figure \ref{figure:2}).
Grâce au lemme \ref{lemma:12.(4.3)}, $a$ et $b$ existent.

\begin{lemma}
\label{lemma:13.(4.4)}
Le groupe engendré par\/ $a$ et\/ $b$ est isomorphe
au produit en couronne restreint\/ $\mathbb Z\wr\mathbb Z$ \textup;
plus précisément\textup,
$$
\langle a,b\rangle=\langle b\rangle\wr\langle a\rangle,\quad
\langle b\rangle\cong\langle a\rangle\cong\mathbb Z.
$$
\end{lemma}


\begin{proof} 
Notons tout d'abord que
$\supp(a^{-k}ba^k)=]\alpha_k;\alpha_{k+1}[$
pour tout $k\in\mathbb Z$
(voir le lemme \ref{lemma:02.(2.2)}).
En particulier, les supports des applications $a^{-k}ba^k$,
$k\in\mathbb Z$, sont deux-à-deux disjoints,
et $\supp(a)=]0;r[$ n'est égal au support d'aucun
élément du groupe $\langle\,a^{-k}ba^k\mid k\in\mathbb Z\,\rangle$.
Donc
$$
\langle\,a^{-k}ba^k\mid k\in\mathbb Z\,\rangle
=\bigoplus_{k\in\mathbb Z}\langle a^{-k}ba^k\rangle
\quad\text{et}\quad
\langle\,a^{-k}ba^k\mid k\in\mathbb Z\,\rangle\cap\langle a\rangle
=\{\id\},
$$
d'où
$\langle a,b\rangle
=\langle\,a^{-k}ba^{k}\mid k\in\mathbb Z\,\rangle\rtimes\langle a\rangle
=\langle b\rangle\wr\langle a\rangle$.
Comme $\F$ est sans torsion,
$\langle a\rangle\cong\langle b\rangle\cong\mathbb Z$.
\end{proof}

Dans le reste de cette section,
$G$ est un sous-groupe de $\V(r,\mathbb R^*_+,\mathbb R)$
tel que
$$
G\cap\F(r,\mathbb R^*_+,\mathbb R)=\F=\F(r,\Lambda,A)
$$
(il est même possible de généraliser certains résultats de cette
section au cas où $G$ satisferait une condition plus faible que celle-ci).

\begin{proposition}
\label{proposition:02.(4.5)}
Soient\/ $a$\textup, $b$\textup, et\/ $G$
les éléments et le groupe définis ci-dessus\textup.
Alors il existe une formule\/ $\phi$ du premier ordre
avec\/ $a$ et\/ $b$ comme paramètres et avec une seule variable libre
tel que\/ $\langle a,b\rangle$
soit défini dans\/ $G$ par\/ $\phi$\textup.
De plus\textup, on peut choisir $\phi$
en fonction seulement de\/ $a$\textup, $b$\textup, et\/ $\F$
\textup(sans avoir connaître\/ $G$ tout entier\/\textup{).}
\end{proposition}

Afin de trouver une telle $\phi$ et ainsi prouver cette proposition,
choisissons tout d'abord $c,d\in\F$ et
$s,t\in\mathbb N$ tels que
\begin{enumerate}
\item
	$c$ soit un générateur du centralisateur de $a$
	dans $\F$,
\item
	$d$ soit un générateur du centralisateur de $b$
	dans $\F_{]\alpha_0;\alpha_1[}$,
\item
	$a=c^s$ et $b=d^t$.
\end{enumerate}
Les éléments $c$ et $d$ et les nombres $s$ et $t$ existent d'après
le lemme \ref{lemma:11.(4.2)}.
D'après le lemme \ref{lemma:04.(2.4)}, $\supp(c)=]0;r[$
et $\supp(d)=]\alpha_0;\alpha_1[$.
Clairement $c,d\in\F^\uparrow$.

Notons que $a^{-k}da^k\in\F^\uparrow$ et que
$\supp(a^{-k}da^k)=]\alpha_k;\alpha_{k+1}[$ pour tout $k$.
Comme pour tout $k$ la conjugaison par $a^k$ est un automorphisme
de $\F$ qui envoie
$\F_{]\alpha_0;\alpha_1[}$ sur
$\F_{]\alpha_{k};\alpha_{k+1}[}$,
$a^{-k}da^k$ est un générateur du centralisateur de $a^{-k}ba^k$
dans $\F_{]\alpha_{k};\alpha_{k+1}[}$
pour tout~$k$.

\begin{lemma}
\label{lemma:14.(4.6)}
Le centralisateur de\/ $\{\,a^{-k}ba^{k}\mid k\in\mathbb Z\,\}$
dans\/ $\F$ est engendré par\/
$\{\,a^{-k}da^{k}\mid k\in\mathbb Z\,\}$\textup.
\end{lemma}
\begin{proof}
L'inclusion
$$
C_{\F}(\{\,a^{-k}ba^{k}\mid k\in\mathbb Z\,\})
\supset\langle\,a^{-k}da^{k}\mid k\in\mathbb Z\,\rangle
$$
est évidente, il reste à montrer l'inclusion inverse.
Alors soit $x$ un élément arbitraire de $\F$
qui commute avec tout $a^{-k}ba^{k}$, $k\in\mathbb Z$.

D'après le lemme \ref{lemma:02.(2.2)}, $x$ permute chacun des intervalles
$]\alpha_k;\alpha_{k+1}[$, $k\in\mathbb Z$.
Par continuité et monotonicité, $(\alpha_k)x=\alpha_k$ pour tout $k$.
Pour tout $k\in\mathbb Z$, soit $y_k$ la permutation de $[0;r[$
telle que
$$
y_k|_{[\alpha_k;\alpha_{k+1}]}=x|_{[\alpha_k;\alpha_{k+1}]}
\quad\text{et}\quad
y_k|_{[0;\alpha_k]\cup[\alpha_{k+1};r[}
=\id_{[0;\alpha_k]\cup[\alpha_{k+1};r[}.
$$
Alors pour tout $k\in\mathbb Z$,
$y_k\in\F_{]\alpha_{k};\alpha_{k+1}[}$ et
$y_k$ commute avec $a^{-k}ba^{k}$.
Donc $y_k\in\langle a^{-k}da^{k}\rangle$ pour tout~$k$.

Choisissons $\beta,\gamma\in]0;r[$ tels que
$x$ soit affine sur $[0;\beta]$ et sur $[\gamma;r[$.
Alors $x|_{[0;\beta]\cup[\gamma;r[}=\id_{[0;\beta]\cup[\gamma;r[}$.
Choisissons $n\in\mathbb N$ tel que $\alpha_{-n}\in]0;\beta]$ et
$\alpha_{n+1}\in[\gamma;r[$.
Alors $\supp(x)\subset]\alpha_{-n};\alpha_{n+1}[$, et donc
$$
x=y_{-n}y_{-n+1}\dotsm y_{n-1}y_{n}
\in\langle\,a^{-k}da^{k}\mid k\in\mathbb Z\,\rangle.
$$
\end{proof}

\begin{lemma}
\label{lemma:15.(4.7)}
Le centralisateur de l'élément\/ $a$ et le centralisateur de
l'ensemble\/ $\{\,a^{-k}ba^{k}\mid k\in\mathbb Z\,\}$
sont conservés en passant
de\/ $\F$ à\/~$G$\textup.
\end{lemma}
\begin{proof}
Tous éléments de $C_{G}(a)$ et de
$C_{G}(\{\,a^{-k}ba^{k}\mid k\in\mathbb Z\,\})$ sont continus
d'après les lemmes \ref{lemma:08.(3.1)} et \ref{lemma:09.(3.2)}.
Comme un élément de $G$ appartient à $\F$
si et seulement si il est continu, la preuve est terminée.
\end{proof}

Le groupe $\langle c\rangle$ est définissable dans $G$
avec le paramètre $a$
puisqu'il est le centralisateur de $a$
(voir le lemme \ref{lemma:15.(4.7)}).
Le groupe $\langle a\rangle$ est définissable dans $G$
avec le même paramètre parce que
$$
\langle a\rangle=\bigl\{\,x^s\bigm|x\in\langle c\rangle\,\bigr\}.
$$
Donc l'ensemble $\{\,a^{-k}ba^{k}\mid k\in\mathbb Z\,\}$
est définissable dans $G$ avec les paramètres
$a$ et $b$, et ainsi son centralisateur l'est.
Par les lemmes \ref{lemma:14.(4.6)} et \ref{lemma:15.(4.7)},
le centralisateur de l'ensemble $\{\,a^{-k}ba^{k}\mid k\in\mathbb Z\,\}$
dans $G$ est le groupe
$\langle\,a^{-k}da^{k}\mid k\in\mathbb Z\,\rangle$.
Comme
$$
\langle\,a^{-k}ba^{k}\mid k\in\mathbb Z\,\rangle
=\bigl\{\,x^t\bigm|
x\in\langle\,a^{-k}da^{k}\mid k\in\mathbb Z\,\rangle\,\bigr\},
$$
le groupe $\langle\,a^{-k}ba^{k}\mid k\in\mathbb Z\,\rangle$ est
définissable avec les mêmes paramètres.
Le groupe $\langle a,b\rangle$ est définissable avec
les paramètres $a$ et $b$
puisqu'il est le produit semi-direct de $\langle a\rangle$
et $\langle\,a^{-k}ba^{k}\mid k\in\mathbb Z\,\rangle$.

La formule suivante définit $\langle a,b\rangle$ dans $G$
et ne dépend que de $a$, $b$, $s$, et~$t$ :
\begin{align*}
\phi(x)=\ulcorner\Bigl(\exists y,z\Bigr)
\Bigl(x=y^{s}z^{t}&\wedge ya=ay\\
&\wedge\Bigl(\forall w\Bigr)
\Bigl(wa=aw\rightarrow zw^{-s}bw^{s}=w^{-s}bw^{s}z\Bigr)\Bigr)\urcorner.
\end{align*}

Nous avons démontré la proposition \ref{proposition:02.(4.5)}.
Nous en déduisons le corollaire suivant.

\begin{corollary}
\label{corollary:01.(4.8)-proposition:02.(4.5)}
Le groupe\/ $\F$ a des sous-groupes isomorphes
à\/ $\mathbb Z\wr\mathbb Z$
et définissables avec paramètres dans\/ $\F$\textup,
dans\/ $\T$\textup, et dans\/~$\V$\textup.
\end{corollary}



\section{Interprétations de l'Arithmétique}
\label{section:interpretationarithmetic}

Dans cette section, nous complétons notre démonstration
de l'interprétabilité de l'Arithmétique
dans $\F$, $\T$, et $\V$ avec paramètres.
En outre, au cas du groupe $F$, ou plus généralement
du groupe $\F$ avec $\Lambda$ cyclique non trivial,
nous présentons une interprétation de l'Arithmétique
qui ne nécessite pas de paramètres.

Apparemment il est bien connu aux spécialistes
que chaque groupe virtuellement résoluble de type fini qui n'est
pas virtuellement abélien interprète l'Arithmétique
(voir \cite{Noskov:1983:etkpprg-rus,Noskov:1984:etfgasg-eng},
et aussi \cite{DelonSimonetta:1998:uwpspsf}).
Pour confort du lecteur, nous presentons ici notre preuve autonome
pour le groupe $\mathbb Z\wr\mathbb Z$.

\begin{lemma}
\label{lemma:16.(5.1)}
Le groupe\/ $\mathbb Z\wr\mathbb Z$ interprète
l'Arithmétique avec paramètres\textup.
Plus précisément\textup, soit\/
$\mathbb Z\wr\mathbb Z=\langle b\rangle\wr\langle a\rangle$\textup,
$\langle b\rangle\cong\langle a\rangle\cong\mathbb Z$\textup,
alors la bijection\/
$$
f\colon\{\,a^n\mid n\in\mathbb N\,\}\to\mathbb N,\quad a^n\mapsto n
$$
est une interprétation de\/
$(\mathbb N,+,\times)$ dans\/ $(\mathbb Z\wr\mathbb Z,\times)$
avec paramètres\textup.
\end{lemma}
\begin{proof}
Notons
$$
G=\mathbb Z\wr\mathbb Z=\langle b\rangle\wr\langle a\rangle
\quad\text{et}\quad
H=\langle\,a^{-k}ba^{k}\mid k\in\mathbb Z\,\rangle.
$$
Rappelons que
$G=H\rtimes\langle a\rangle$,
et que $H$ est un groupe abélien libre avec la base
$(a^{-k}ba^{k})_{k\in\mathbb Z}$.
On peut vérifier facilement que $\langle a\rangle=C_G(a)$,
et que $H=C_G(b)$.

Considérons la bijection
$$
g\colon \langle a\rangle\to\mathbb Z,\quad a^n\mapsto n.
$$
Il suffira de démontrer que $g$ est une interprétation de
$(\mathbb Z,+,\times)$ dans $(G,\times)$.
En effet, $(\mathbb N,+,\times)$ est une sous-structure de
$(\mathbb Z,+,\times)$, et
$\mathbb N$ est $0$-définissable
dans $(\mathbb Z,+,\times)$ parce que,
d'après le théorème des quatre carrés de Lagrange,
tout entier positif est la somme de quatre carrés.

L'ensemble de départ de $g$ est le centralisateur de $a$,
donc définissable.
L'opération induite sur $\langle a\rangle$
par l'addition de $\mathbb Z$ via $g$ est
tout simplement la restriction de la multiplication de $G$,
donc $0$-définissable.
Il reste à démontrer que l'opération induite sur $\langle a\rangle$
par la multiplication de $\mathbb Z$ via $g$ est définissable.

Remarquons les faits suivants :
\begin{enumerate}
\item
	pour tout $x\in G\setminus H$, $C_G(x)$ est cyclique,
\item
	pour tout $n\in\mathbb Z\setminus\{0\}$,
	$C_G(ba^n)=\langle ba^n\rangle$,
\item
	pour tout $n\in\mathbb Z$,
	$HC_G(ba^n)=H\langle a^n\rangle$.
\end{enumerate}
(Le deuxième fait est dû à l'homomorphisme
$G\to\mathbb Z,\ a\mapsto0,\ b\mapsto1$.)

Notons $|$ la relation de la divisibilité
dans $\mathbb Z$.
Observons que pour tous $m,n\in\mathbb Z$,
$$
m|n
\:\Leftrightarrow\:
HC_G(ba^m)\supset HC_G(ba^n).
$$
La relation $HC_G(bx)\supset HC_G(by)$ entre $x,y\in G$
s'exprime par une formule du premier ordre
avec le paramètre $b$.
Donc la relation induite sur $\langle a\rangle$ par $|$ via $g$
est définissable.

La multiplication dans $\mathbb Z$ est définissable
à partir de l'addition, la divisibilité, et la constante $1$
grâce aux équivalences suivantes satisfaites dans~$\mathbb Z$ :
\begin{align*}
n=k(k+1)&\leftrightarrow
\Bigl(\forall m\Bigr)\Bigl(n|m\leftrightarrow
k|m\wedge(k+1)|m\Bigr)\wedge(2k+1)|(2n-k),\\
n=kl&\leftrightarrow
(k+l)(k+l+1)=k(k+1)+l(l+1)+2n
\end{align*}
(voir \cite[\S5a]{Robinson:1951:ur} pour les détails).
Donc l'opération induite sur $\langle a\rangle$
par la multiplication de $\mathbb Z$
via $g$ est définissable.
\end{proof}

Le théorème A (voir l'introduction)
est un corollaire de la proposition \ref{proposition:02.(4.5)}
et du lemme \ref{lemma:16.(5.1)}.
Afin de démontrer le théorème B, nous construirons
des $0$-interprétations nouvelles de l'Arithmétiques dans
des groupes de genre $\F$.


\begin{proposition}
\label{proposition:03.(5.2)}
Si\/ $\Lambda$ est cyclique\textup,
$\Lambda=\langle p\rangle$\textup,
alors l'application\/
$$
f\colon \{\,x\in\F\mid(0)x^{\prime+}=(r)x^{\prime-}>1\,\}\to\mathbb N,
\quad x\mapsto\log_p((0)x^{\prime+})
$$
est une interprétation de\/
$(\mathbb N,+,\times)$ dans\/ $(\F,\times)$
sans paramètres\textup.
\end{proposition}

Une des idées principales de la preuve de la proposition
\textup{\ref{proposition:03.(5.2)}} est l'usage des
centralisateurs de paires d'éléments.%
\footnote{Les centralisateurs dans $\F(r,\mathbb R^*_+,\mathbb R)$ ont été
décrits par Brin et Squier \cite{BrinSquier:2001:pcrcgplhrl}.
Collin Bleak et autres \cite{BleakGGHMNS:pp2007:} ont
récemment annoncé une classification de tous
les centralisateurs dans
$\T(1,\langle n\rangle,\mathbb Z[\frac{1}{n}])$
et $\V(1,\langle n\rangle,\mathbb Z[\frac{1}{n}])$, $n=2,3,\dotsc$.}
Le lemme suivant est similaire au théorème 5.5
dans~\cite{BrinSquier:2001:pcrcgplhrl}.

\begin{lemma}
\label{lemma:17.(5.3)}
Soit\/ $H$ un sous-groupe de\/ $\F$\textup.
Alors\/ $H$ est le centralisateur d'un élément
si et seulement si\/ $H$ se décompose en un produit direct
des sous-groupes\/ $H_1,\dotsc,H_n$\textup,
$n\in\mathbb N$\textup, tels qu'il existe
$\alpha_0,\dotsc,\alpha_n\in A$ tels que\/ \textup:
\begin{enumerate}
\item
	$0=\alpha_0<\alpha_1<\dotsb<\alpha_n=r$ \textup;
\item
	pour tout\/ $i=1,\dotsc,n$\textup,
	soit\/ $H_i=\F_{]\alpha_{i-1};\alpha_i[}$\textup,
	soit il existe\/ $x$ dans\/ $\F_{]\alpha_{i-1};\alpha_i[}$ tel que\/
	$H_i$ soit le centralisateur de $x$ dans\/
	$\F_{]\alpha_{i-1};\alpha_i[}$\textup, et que\/
	$H_i=\langle x\rangle$ \textup;
\item
	pour tout\/ $i=1,\dotsc,n-1$\textup,
	si\/ $H_i=\F_{]\alpha_{i-1};\alpha_i[}$\textup, alors\/
	$H_{i+1}\ne\F_{]\alpha_i;\alpha_{i+1}[}$\textup.
\end{enumerate}
\end{lemma}
\begin{proof}
Soient $x\in\F$ et $H=C_\F(x)$.
Choisissons $\alpha_0,\dotsc,\alpha_n$
tels que $0=\alpha_0<\alpha_1<\dotsb<\alpha_n=r$ et que
$$
\{\alpha_1,\dotsc,\alpha_{n-1}\}
=\{\,\alpha\in]0;r[\cap A\cap\fix(x)\mid
(\alpha)x^{\prime-}\ne1\text{ ou }(\alpha)x^{\prime+}\ne1\,\}.
$$
Pour tout $i=1,\dotsc,n$,
soit $x_i\in\F_{]\alpha_{i-1};\alpha_i[}$ tel que
$x_i|_{]\alpha_{i-1};\alpha_i[}=x|_{]\alpha_{i-1};\alpha_i[}$,
et soit $H_i$ le centralisateur de $x_i$
dans $\F_{]\alpha_{i-1};\alpha_i[}$.
Notons que $x=x_1\dotsm x_n$.
Pour tout $i=1,\dotsc,n$, si $x_i\ne\id$, alors $H_i$ est cyclique
(voir le lemme \ref{lemma:11.(4.2)}).

Pareillement au lemme \ref{lemma:03.(2.3)}, il est facile à prouver que
tout élément $y$ de $H$ permute
l'ensemble $\{\alpha_0,\dotsc,\linebreak[0]\alpha_{n-1}\}$,
et donc, comme cet ensemble est fini, $y$ fixe tous ses éléments.
Donc $H=H_1\times\dotsm\times H_n$.

Réciproquement, supposons que $H=H_1\times\dotsm\times H_n$,
$n\in\mathbb N$, et que les sous-groupes $H_1,\dotsc,H_n$
et les points $\alpha_0,\dotsc,\alpha_n\in A$
soient comme dans l'énoncé de ce lemme.
Pour tout $i=1,\dotsc,n$, soit $x_i\in\F_{]\alpha_{i-1};\alpha_i[}$
tel que $H_i$ soit le centralisateur de $x_i$
dans $\F_{]\alpha_{i-1};\alpha_i[}$.
Alors $H=C_\F(x_1\dotsm x_n)$.
\end{proof}

\begin{lemma}
\label{lemma:18.(5.4)}
Soit\/ $H$ un sous-groupe de\/ $\F$\textup.
Alors\/ $H$ est le centralisateur d'une paire d'éléments\/
\textup(éventuellement identiques\/\textup)
si et seulement si\/ $H$ se décompose en un produit direct
des sous-groupes\/ $H_1,\dotsc,H_n$\textup,
$n\in\mathbb N$\textup, tels qu'il existe
$\alpha_0,\dotsc,\alpha_n\in A$ tels que\/ \textup:
\begin{enumerate}
\item
	$0=\alpha_0<\alpha_1<\dotsb<\alpha_n=r$ \textup;
\item
	pour tout\/ $i=1,\dotsc,n$\textup,
	soit\/ $H_i=\{\id\}$\textup,
	soit\/ $H_i=\F_{]\alpha_{i-1};\alpha_i[}$\textup,
	soit il existe\/ $x$ dans\/ $\F_{]\alpha_{i-1};\alpha_i[}$ tel que\/
	$H_i$ soit le centralisateur de $x$ dans\/
	$\F_{]\alpha_{i-1};\alpha_i[}$\textup, et que\/
	$H_i=\langle x\rangle$ \textup;
\item
	pour tout\/ $i=1,\dotsc,n-1$\textup,
	si\/ $H_i=\F_{]\alpha_{i-1};\alpha_i[}$\textup, alors\/
	$H_{i+1}\ne\F_{]\alpha_i;\alpha_{i+1}[}$\textup.
\end{enumerate}
\end{lemma}
\begin{proof}
Ce lemme est un corollaire facile du lemme \ref{lemma:17.(5.3)}
et du fait que pour tous $\alpha,\beta\in A$
tels que $0<\alpha<\beta<r$, il existe $x,y\in\F$ tels que
$\supp(x)=\supp(y)=]\alpha;\beta[$ et que $xy\ne yx$,
et donc le centralisateur de $\{x,y\}$ dans $\F_{]\alpha;\beta[}$
est trivial (voir les lemmes \ref{lemma:11.(4.2)} et \ref{lemma:12.(4.3)}).
\end{proof}

\begin{proof}
[Démonstration de la proposition\/ \textup{\ref{proposition:03.(5.2)}}]
Sans perte de généralité, supposons que $p>1$.

Notons $\F^\circ$ le sous-groupe de $\F$ formé des éléments qui
sont l'identité aux voisinages de $0$ et de $r$ :
$$
\F^\circ=\{\,x\in\F\mid
(0)x^{\prime+}=(r)x^{\prime-}=1\,\}.
$$
Notons $B$ l'ensemble de départ de $f$ :
$$
B=\{\,x\in\F\mid(0)x^{\prime+}=(r)x^{\prime-}>1\,\}.
$$
Remarquons que pour deux éléments $x$ et $y$ de $B$,
$f(x)=f(y)$ si et seulement si $xy^{-1}\in\F^\circ$
(voir le lemme~\ref{lemma:10.(4.1)}).

Le lemme \ref{lemma:12.(4.3)} permet de conclure que
$f$ est une surjection sur $\mathbb N$
(voir la figure \ref{figure:3} par exemple).
\begin{figure}
\includegraphics{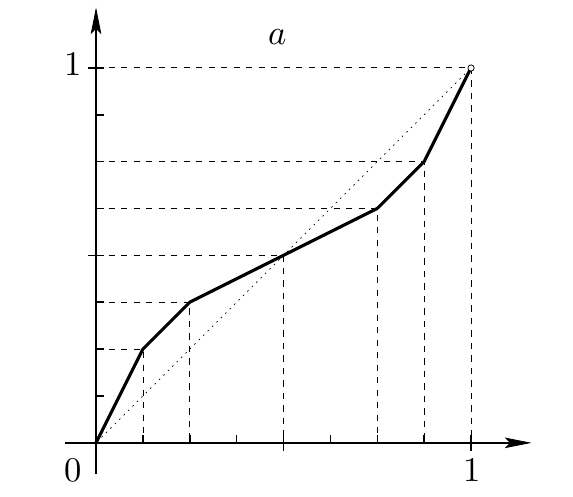}
\caption{Un élément $a$ de $F$ tel que $(0)a^{\prime+}=(1)a^{\prime-}=2$.}
\label{figure:3}
\end{figure}
Il nous suffira de démontrer que l'ensemble $B$, le groupe $\F^\circ$,
et les relations induites sur $B$
via $f$ par l'addition et la divisibilité de $\mathbb N$
sont tous $0$-définissables (voir la remarque \ref{remark:04.(2.10)}).
En effet, la multiplication est $0$-définissable dans $\mathbb N$
à partir de l'addition et de la divisibilité ---
voir \cite[\S4b]{Robinson:1951:ur} ou la preuve du
lemme~\ref{lemma:16.(5.1)}.

Notons
\begin{align*}
S&=\{\,\F_{]\alpha;\beta[}\mid\alpha,\beta\in A,\
0\le\alpha<\beta\le r\,\},\\
S_0&=\{\,\F_{]\alpha;\beta[}\mid\alpha,\beta\in A,\
0<\alpha<\beta<r\,\},\\
S_1&=\{\,\F_{]0;\beta[}\mid\beta\in A,\ 0<\beta<r\,\}\\
&\qquad\qquad\cup
\{\,\F_{]\alpha;r[}\mid\alpha\in A,\ 0<\alpha<r\,\}.
\end{align*}
Alors non seulement la famille $S$ est \emph{uniformément définissable},
mais surtout
il existe deux formules du premier ordre
$\phi(x_1,x_2,x_3)$ et $\psi(x_1,x_2)$
(dans le langage des groupes) \emph{sans paramètres\/}
telles que :
\begin{enumerate}
\item
	pour tous $\alpha,\beta\in A$ tels que
	$0\le\alpha<\beta\le r$, il existe $x,y\in\F$ tels que
	$\F\models\psi(x,y)$ et que
	$\F_{]\alpha;\beta[}=\{\,z\in\F\mid\F\models\phi(x,y,z)\,\}$ ;
\item
	pour tous $x,y\in\F$ tels que
	$\F\models\psi(x,y)$, il existe $\alpha,\beta\in A$ tels que
	$0\le\alpha<\beta\le r$ et que
	$\F_{]\alpha;\beta[}=\{\,z\in\F\mid\F\models\phi(x,y,z)\,\}$.
\end{enumerate}
En effet, une partie de $\F$ est de la forme
$\F_{]\alpha;\beta[}$, où $\alpha,\beta\in A$ et
$0\le\alpha<\beta\le r$,
si et seulement si elle est le centralisateur d'une paire,
qui n'est pas abélien, et qui
ne se décompose pas comme le produit direct de deux autres
centralisateurs de paires (voir le lemme \ref{lemma:18.(5.4)}).
Tout cela s'exprime au premier ordre.
Les familles $S_0$ et $S_1$ de parties de $\F$ sont
«~uniformément définissables sans paramètres~» au même sens que $S$
parce que
\begin{align*}
S_0&=\{\,H_1\cap H_2\mid
H_1,H_2\in S,\ H_1\not\subset H_2,\text{ et }H_2\not\subset H_1\,\},\\
S_1&=S\setminus(S_0\cup\{\F\}).
\end{align*}
Comme $\F^\circ=\bigcup S_0$, ce sous-groupe est $0$-définissable.%
\footnote{Dans le cas où $\F=F$,
il existe une preuve plus naturelle pour montrer que
$\F^\circ$ est définissable dans $\F$,
comme dans ce cas $\F^\circ=F^\circ=[F,F]=[\F,\F]$
(voir \cite[théorème 4.1]{CannonFloydParry:1996:inRTg}),
et que tout élément de $[\F,\F]$ est le produit de deux commutateurs
(voir l'appendice).
En général, $\F^\circ$ et $[\F,\F]$ ne sont pas égaux
(voir \cite[section 4D]{Brown:1987:fpg}).}

Notons
\begin{align*}
E&=\{\,x\in\F\mid(0)x^{\prime+}\ne1\text{ ou }(r)x^{\prime-}\ne1\,\}
=\F\setminus\bigcup S_0=\F\setminus\F^\circ,\\
E_2&=\{\,x\in\F\mid(0)x^{\prime+}\ne1\text{ et }(r)x^{\prime-}\ne1\,\}
=\F\setminus\bigcup S_1=E\setminus\bigcup S_1.
\end{align*}
Ces ensembles sont $0$-définissables.

Notons
\begin{align*}
P^+&=\{\,x\in\F\mid(0)x^{\prime+}>1\text{ et }(r)x^{\prime-}>1\,\},\\
P^-&=\{\,x\in\F\mid(0)x^{\prime+}<1\text{ et }(r)x^{\prime-}<1\,\},
\end{align*}
et $P=P^+\cup P^-$.
Ces ensembles sont $0$-définissables :
pour tout $x\in\F$,
$$
x\in P^+
\:\Leftrightarrow\:
\Bigl(\exists X\in S_0\Bigr)\,\Bigl(\forall Y\in S_1\Bigr)\,
\Bigl(Y\supset X\rightarrow x^{-1}Yx\subsetneqq Y\Bigr),
$$
et
$$
x\in P^-
\:\Leftrightarrow\:
\Bigl(\exists X\in S_0\Bigr)\,\Bigl(\forall Y\in S_1\Bigr)\,
\Bigl(Y\supset X\rightarrow x^{-1}Yx\supsetneqq Y\Bigr).
$$
Il a été utilisé ici que pour tous $\alpha,\beta$
et pour tout $x\in\F$,
$$
x^{-1}\F_{]\alpha;\beta[}x=\F_{](\alpha)x;(\beta)x[}.
$$

Notons
$$
U=\{\,x\in\F\mid(0)x^{\prime+}=((r)x^{\prime-})^{-1}\in\{p^{\pm1}\}\,\}.
$$
Alors $U$ est $0$-définissable : $U\subset E_2\setminus P$, et
pour tout $x\in E_2\setminus P$,
\begin{multline*}
x\in U
\:\Leftrightarrow\:
\Bigl(\forall y\in P\Bigr)\,
\Bigl(\forall z\in\F^\circ\Bigr)\,
\Bigl(\exists w_1,w_2\in C_\F(xz)\Bigr)\\
\Bigl((w_1w_2^{-1}\in E_2)\wedge (yw_1,yw_2\notin E_2)\Bigr).
\end{multline*}
Afin de faciliter la lecture de la dernière formule, remarquons que
$xy^{-1}\notin E_2$ signifie exactement que soit
$(0)x^{\prime+}=(0)y^{\prime+}$, soit $(r)x^{\prime-}=(r)y^{\prime-}$.
Alors l'implication «~$\Rightarrow$~» est facile à démontrer :
on peut toujours trouver de tels $w_1,w_2$ même
dans $\langle xz\rangle$.
La direction «~$\Leftarrow$~» est moins évidente, démontrons-la par
contraposé.

Soit $x\in(E_2\setminus P)\setminus U$.
Sans perte de généralité, supposons que
$(0)x^{\prime+}>p$, et donc $(r)x^{\prime-}<1$.
Choisissons $y\in P$ tel que $(0)y^{\prime+}=p$
(voir le lemme \ref{lemma:12.(4.3)}).
Soit $\gamma\in]0;r[\cap A$, choisissons $x_1,x_2\in\F$ tels que :
\begin{enumerate}
\item
	$\supp(x_1)=]0;\gamma[$, $\supp(x_2)=]\gamma;r[$,
\item
	$(0)x_1^{\prime+}=(0)x^{\prime+}$,
	$(\gamma)x_1^{\prime-}=p^{-1}$,
\item
	$(r)x_2^{\prime-}=(r)x^{\prime-}$,
	$(\gamma)x_2^{\prime+}=p$.
\end{enumerate}
Alors le centralisateur de $x_1x_2$ est le produit direct
$\langle x_1\rangle\times\langle x_2\rangle$
(voir les lemmes \ref{lemma:03.(2.3)} et \ref{lemma:10.(4.1)}).
Choisissons $z\in\F^\circ$ tel que $xz=x_1x_2$.
Supposons maintenant que
$$
w_1,w_2\in C_\F(xz)=\langle x_1\rangle\times\langle x_2\rangle.
$$
Alors $(0)(yw_1)^{\prime+}\ne1$ et $(0)(yw_2)^{\prime+}\ne1$.
Supposons que $yw_1,yw_2\notin E_2$.
Alors $(r)(yw_1)^{\prime-}=(r)(yw_2)^{\prime-}=1$,
et donc $(r)(w_1w_2^{-1})^{\prime-}=1$,
soit encore $w_1w_2^{-1}\notin E_2$.
Donc $U$ est en effet $0$-définissable.

Pour tout $x\in P$,
$$
(0)x^{\prime+}=(r)x^{\prime-}
\:\Leftrightarrow\:
\Bigl(\forall y\in U\Bigr)\,
\Bigl(\exists z\in C_\F(y)\Bigr)\,
\Bigl(xz,xz^{-1}\notin E_2\Bigr).
$$
Donc, l'ensemble $B$ est $0$-définissable.

La $f$-préimage du graphe de l'addition de $\mathbb N$
est $0$-définissable comme elle est le graphe de la
multiplication modulo $\F^\circ$ : pour tous $x,y,z\in B$,
$$
(x)f+(y)f=(z)f
\:\Leftrightarrow\:
xyz^{-1}\in\F^\circ.
$$

Il nous reste à montrer que la $f$-préimage du graphe de
la divisibilité de $\mathbb N$ est $0$-définissable.
C'est en effet le cas :
pour tous $x,y\in B$,
$$
(x)f|(y)f
\:\Leftrightarrow\:
\Bigl(\forall z\in\F^\circ\Bigr)\,
\Bigl(\exists w\in C_\F(xz)\Bigr)\,
\Bigl(yw\in\F^\circ\Bigr).
$$
L'implication «~$\Leftarrow$~» ici est la moins évidente.
Afin de la démontrer par contraposé, on peut prendre
$\gamma\in]0;r[\cap A$, choisir $x_1,x_2\in\F$ tels que :
\begin{enumerate}
\item
	$\supp(x_1)=]0;\gamma[$, $\supp(x_2)=]\gamma;r[$,
\item
	$(0)x_1^{\prime+}=(0)x^{\prime+}$,
	$(\gamma)x_1^{\prime-}=p^{\pm1}$,
\item
	$(r)x_2^{\prime-}=(r)x^{\prime-}$,
	$(\gamma)x_2^{\prime+}=p^{\pm1}$,
\end{enumerate}
et choisir $z\in\F^\circ$ tel que $xz=x_1x_2$
(et donc $C_\F(xz)=\langle x_1\rangle\times\langle x_2\rangle$).
\end{proof}

Le théorème B
est une conséquence de la proposition \textup{\ref{proposition:03.(5.2)}}.


\section{Indécidabilité}
\label{section:undecidability}

Dans cette section, nous déduisons de
\cite[théorème 9]{MostowskiRobTar:1971:ueua}
(voir le théorème \ref{theorem:01.(6.3)} ci-dessous)
que la théorie élémentaire de toute structure de signature finie
qui interprète l'Arithmétique avec paramètres
est héréditairement indécidable.

Dans ce qui suit, les constantes sont traitées comme des
fonctions d'\emph{arité\/} $0$,
et également des symboles de constante sont vus comme
un cas particulier de symboles de fonction.
Si $\sigma$ est un symbole de relation, de fonction, ou de constante,
notons $\arity(\sigma)$ l'arité de $\sigma$.

Disons qu'une relation $R$ $n$-aire sur un ensemble $B$ est
\emph{compatible\/} avec une relation d'équivalence $E$ sur $B$ si
$R$ est induite par une relation (de la même arité) sur $B/E$,
c'est-à-dire si l'appartenance à $R$ d'une $n$-uple
$(b_1,\dotsc,b_n)$
ne dépend que des classes de $E$-équivalence de $b_1,\dotsc,b_n$.

Soient $\Sigma$ et $\Gamma$ deux signatures, et soit
$N$ une $\Gamma$-structure.
Soient $n\in\mathbb N$, $B$ une partie définissable de $N^{n}$,
et $E$ une relation d'équivalence sur $B$
aussi définissable dans $N$.
Soient $\phi$, $\psi$, et $\xi_\sigma$ pour tout $\sigma\in\Sigma$
des $\Gamma$-formules avec des paramètres extraits de $N$
telles que :
\begin{enumerate}
\item
	$\phi$ définisse $B$,
\item
	$\psi$ définisse $E$,
\item
	pour tout symbole de relation $\sigma\in\Sigma$,
	la formule $\xi_\sigma$
	définisse une relation sur $B$ compatible avec $E$
	d'arité $\arity(\sigma)n$, et
\item
	pour tout symbole de fonction ou de constante $\sigma\in\Sigma$,
	la formule $\xi_\sigma$
	définisse une relation sur $B$
	d'arité $(\arity(\sigma)+1)n$ qui est compatible avec $E$ et
	qui définit le graphe d'une opération sur $B/E$
	d'arité $\arity(\sigma)$.
\end{enumerate}
Alors notons $\Int{\Sigma}{N,\phi,\psi,(\xi_\sigma)_{\sigma\in\Sigma}}$
la $\Sigma$-structure définie sur $B/E$ par la famille
$(\xi_\sigma)_{\sigma\in\Sigma}$ au sens naturel.

\begin{remark}
\label{remark:07.(6.1)}
Dans la même notation, la projection naturelle
$p\colon B\to B/E$ est une interprétation
de $\Int{\Sigma}{N,\phi,\psi,(\xi_\sigma)_{\sigma\in\Sigma}}$ dans~$N$.
\end{remark}

\begin{lemma}
\label{lemma:19.(6.2)}
Soient\/ $M$ et\/ $N$ deux structures de signatures finies telles
que\/ $\Th(M)$ soit héréditairement indécidable\textup, et
que\/ $N$ interprète\/ $M$ avec paramètres\textup.
Alors\/ $\Th(N)$ est héréditairement indécidable aussi\textup.
\end{lemma}
\begin{proof}
Notons $\Sigma$ la signature de $M$ et $\Gamma$ la signature de~$N$.

Il suffit de considérer le cas où $\Sigma$
ne contient pas de symboles de fonction, ni de symboles de constante.
En effet, soit $\Sigma'$ la signature obtenue à partir de $\Sigma$
en remplaçant
tout symbole de fonction $n$-aire $f$ ($n\ge0$) par un symbole
de relation $(n+1)$-aire $f'$.
Pour toute $\Sigma$-structure $M$, notons $M'$ la $\Sigma'$-structure
sur l'ensemble sous-jacent de $M$ où tout nouveau symbole de relation
de $\Sigma'$ est interprété par le graphe de la fonction dans $M$
nommée par le symbole de $\Sigma$ correspondant,
et tous les autres symboles de $\Sigma'$ sont interprétés dans $M'$
exactement comme dans $M$.
Pour toute $\Sigma$-théorie $S$,
notons $S'$ la $\Sigma'$-théorie
de la classe $\{\,M'\mid M\models S\,\}$.
Il est facile de voir que :
\begin{enumerate}
\item
	la (classe-)application $M\mapsto M'$,
	où $M$ est un $\Sigma$-modèle de $S$,
	est une (classe-) bijection entre
	$\Mod_{\Sigma}S$ et $\Mod_{\Sigma'}S'$ pour toute
	$\Sigma$-théorie $S$,
\item
	si $O_{\Sigma}$ dénote la $\Sigma$-théorie minimale, alors
	l'application $S\mapsto S'$
	de $\Sigma$-théories vers $\Sigma'$-théories
	est une bijection entre toutes les $\Sigma$-théories
	et toutes les $\Sigma'$-théories contenant $O_{\Sigma}'$
	($O_{\Sigma}'$ exprime tout
	simplement que les nouveaux symboles de relation de $\Sigma'$
	s'interprètent par des graphes de fonctions),
\item
	pour toute $\Sigma$-structure $M$,
	un ensemble est $0$-définissable dans $M$ si et seulement
	si il l'est dans $M'$,
	ou autrement dit,
	$\id_M$ est une $0$-interprétation
	de $M$ dans $M'$ et de $M'$ dans~$M$.
\end{enumerate}
Il est facile de fournir un algorithme
qui convertit tout $\Sigma$-énoncé $\phi$ en un $\Sigma'$-énoncé
$\psi$ tel que pour toute $\Sigma$-structure $M$,
$M\models\phi$ si et seulement si $M'\models\psi$,
et il est également facile de fournir un algorithme
qui convertit tout $\Sigma'$-énoncé en
un $\Sigma$-énoncé équivalent dans le même sens.
En conséquence, pour toute $\Sigma$-théorie $S$,
$S'$ est essentiellement indécidable si
et seulement si $S$ l'est.
Donc nous supposons sans perte de généralité que
$\Sigma$ ne contienne que des symboles de relation.

Soit $(n,f)$ une interprétation de $M$ dans $N$.
Soit $a_1,\dotsc,a_m$ une suite des paramètres extraits de $N$
suffisante pour définir l'ensemble de départ et le noyau de $f$ et
la $f$-préimage du graphe de toute relation de $M$ ($\Sigma$ est finie) ;
notons $\bar a=(a_1,\dotsc,a_m)$.
Soient $x_1,\dotsc,x_m$, $y_1,\dotsc,y_n$,
$y_{11},\dotsc,y_{1n}$, $y_{21},\dotsc,y_{2n}$, \dots\
des variables différentes, et notons
$\bar x=(x_1,\dotsc,x_m)$,
$\bar y=(y_1,\dotsc,y_n)$,
$\bar y_1=(y_{11},\dotsc,y_{1n})$, et ainsi de suite.
Soient $\phi=\phi(\bar x,\bar y)$,
$\psi=\psi(\bar x,\bar y_1,\bar y_2)$,
et $\xi_\sigma=\xi_\sigma(\bar x,\bar y_1,\dotsc,\bar y_k)$
pour tout $\sigma\in\Sigma$ d'arité $k$
des $\Gamma$-formules telles que :
\begin{enumerate}
\item
	$\phi(\bar a,\bar y)$ définisse l'ensemble de départ de $f$
	(qui est une partie de $N^n$),
\item
	$\psi(\bar a,\bar y_1,\bar y_2)$ définisse le noyau de $f$
	(qui est une relation d'équivalence sur l'ensemble de départ),
\item
	pour tout symbole $\sigma\in\Sigma$,
	la formule $\xi_\sigma(\bar a,\bar y_1,\dotsc,\bar y_{\arity(\sigma)})$
	définisse la $f$-préimage du graphe de la relation de $M$ nommée
	par~$\sigma$.
\end{enumerate}
Remarquons que la bijection induite par $f$
entre le quotient de son ensemble de départ par son noyau et
son image est un isomorphisme
$$
\Int{\Sigma}{N,\phi(\bar a,\bar y),
\psi(\bar a,\bar y_1,\bar y_2),
(\xi_\sigma(\bar a,\bar y_1,\dotsc,
\bar y_{\arity(\sigma)}))_{\sigma\in\Sigma}}\stackrel{\cong}{\to}M.
$$

Dans ce qui suit, soit $\bar c=(c_1,\dotsc,c_m)$ une
suite de nouveaux symboles de constante.
Écrivons $(\Gamma,\bar c)$ pour noter la signature obtenue
à partir de $\Gamma$ en rajoutant $c_1,\dotsc,c_m$
(comme des symboles de constante).

Soit $\tau=\tau(\bar x)$ une $\Gamma$-formule telle que
la $(\Gamma,\bar c )$-énoncé $\tau(\bar c)$ exprime que :
\begin{enumerate}
\item
	$\psi(\bar c,\bar y_1,\bar y_2)$
	définit une relation d'équivalence sur
	l'ensemble défini par $\phi(\bar c,\bar y)$,
\item
	pour tout symbole de relation $\sigma\in\Sigma$,
	la formule
	$\xi_\sigma(\bar c,\bar y_1,\dotsc,\bar y_{\arity(\sigma)})$
	définit une relation
	sur l'ensemble défini par $\phi(\bar c,\bar y)$
	qui est compatible avec la relation d'équivalence
	définie par $\psi(\bar c,\bar y_1,\bar y_2)$.
\end{enumerate}
Clairement tout cela s'exprime au premier ordre.%
\footnote{On peut prendre comme $\tau(\bar c)$
la conjonction des \emph{conditions d'admissibilité\/} au sens de
\cite[section 5.3]{Hodges:1993:mt}.}
Notons que
$$
N\models\tau(\bar a).
$$


Choisissons une application $t$ récursive
(\ie\ calculable par un algorithme) de l'ensemble des $\Sigma$-énoncés
vers l'ensemble des $\Gamma$-formules avec toutes ses variables libres
comprises parmi $x_1,\dotsc,x_m$ telle que
pour tout $\Sigma$-énoncé $\alpha$ et
toute $(\Gamma,\bar c)$-structure $L$ telle que $L\models\tau(\bar c)$,
$$
\Bigl(L\models\alpha^t(\bar c)\Bigr)
\Leftrightarrow
\Bigl(\Int{\Sigma}{L,\phi(\bar c,\bar y),\psi(\bar c,\bar y_1,\bar y_2),
(\xi_\sigma(\bar c,\bar y_1,\dotsc,
\bar y_{\arity(\sigma)}))_{\sigma\in\Sigma}}\models\alpha\Bigr).
$$
Il est facile de construire une telle $t$ qui utilise la formule 
$\phi$ pour relativiser les quanteurs,
la formule $\psi$ pour remplacer $\ulcorner=\urcorner$, et la formule
$\xi_\sigma$ pour remplacer chaque $\sigma\in\Sigma$.%
\footnote{Dans \cite[section 5.3]{Hodges:1993:mt}
une telle $t$ est dite une \emph{application de réduction}.}
Voici un exemple, où $\sigma$ est un symbole de relation binaire :
\begin{align*}
\alpha&=\ulcorner\Bigl(\forall y_1,y_2,y_3\Bigr)\,
\Bigl(\sigma(y_1,y_2)\wedge\sigma(y_1,y_3)
\rightarrow y_2=y_3\Bigr)\urcorner,\\
\alpha^t(\bar x)&=\ulcorner\Bigl(\forall\bar y_1,\bar y_2,\bar y_3\Bigr)\,
\Bigl(\phi(\bar x,\bar y_1)\wedge\phi(\bar x,\bar y_2)
\wedge\phi(\bar x,\bar y_3)\\
&\qquad\qquad\rightarrow
\Bigl(\xi_\sigma(\bar x,\bar y_1,\bar y_2)
\wedge\xi_\sigma(\bar x,\bar y_1,\bar y_3)
\rightarrow\psi(\bar x,\bar y_2,\bar y_3)\Bigr)\Bigr)\urcorner.
\end{align*}
(Comme à l'accoutumée, nous ne montrons pas toutes les parenthèses ;
elles devraient être ajoutées selon les règles standards.)

Notons les propriétés suivantes de $t$ :
\begin{enumerate}
\item
	pour tout $\Sigma$-énoncé $\alpha$,
	$$
	\Bigl(M\models\alpha\Bigr)\Leftrightarrow
	\Bigl(N\models\alpha^t(\bar a)\Bigr) ;
	$$
\item
	pour tout $\Sigma$-énoncé $\alpha$,
	$$
	\Bigl(\vdash_\Sigma\alpha\Bigr)\Rightarrow
	\Bigl(\tau(\bar c)\vdash_{(\Gamma,\bar c)}\alpha^t(c)\Bigr) ;
	$$
\item
	pour tous $\Sigma$-énoncés $\alpha$ et $\beta$,
	\begin{align*}
	\tau(\bar c)&\vdash_{(\Gamma,\bar c)}(\alpha\wedge\beta)^t(c)
	\leftrightarrow\alpha^t(c)\wedge\beta^t(c),\\
	\tau(\bar c)&\vdash_{(\Gamma,\bar c)}(\neg\alpha)^t(c)
	\leftrightarrow\neg\alpha^t(c),
	\end{align*}
	et de même pour les autres opérations booléennes ;
\item
	pour toute $(\Gamma,\bar c)$-théorie $T$ telle que
	$T\vdash_{(\Gamma,\bar c)}\tau(\bar c)$, l'ensemble
	$$
	\{\,\text{$\Sigma$-énoncé }\alpha
	\mid T\vdash_{(\Gamma,\bar c)}\alpha^t(\bar c)\,\}
	$$
	est une $\Sigma$-théorie.
\end{enumerate}


Supposons maintenant que $\Th(N)$ ne soit pas héréditairement
indécidable.
Alors soit $T$ une $\Gamma$-sous-théorie décidable de $\Th(N)$.
Soit $S$ l'ensemble de toutes les $\Gamma$-formules $\alpha(\bar x)$
telles que
$$
T\vdash_\Gamma\Bigl(\forall\bar x\Bigr)\,
\Bigl(\tau(\bar x)\rightarrow\alpha(\bar x)\Bigr).
$$
Alors $\tau\in S$, $N\models\alpha(\bar a)$ pour toute $\alpha\in S$,
$S$ est un ensemble récursif (décidable),
et $\{\,\alpha(\bar c)\mid\alpha(\bar x)\in S\,\}$
est une $(\Gamma,\bar c)$-théorie.
Soit $U$ la préimage de $S$ sous $t$.
Alors $U\subset\Th(M)$, $U$ est une $\Sigma$-théorie,
et $U$ est décidable contrairement à
l'indécidabilité héréditaire de~$\Th(M)$.
\end{proof}


\begin{theorem}[Mostowski, Tarski, \cite{MostowskiTarski:1949:uaitr}]
\label{theorem:01.(6.3)}
La théorie élémentaire de l'Arithmétique\/ $(\mathbb N,+,\times)$ a une
sous-théorie essentiellement indécidable et finiment axiomatisée\textup.
\end{theorem}

Pour une preuve amélorée de ce fait,
voir \cite[théorème 9]{MostowskiRobTar:1971:ueua}.


\begin{proposition}
\label{proposition:04.(6.4)}
Soit\/ $M$ une structure de signature finie
qui interprète l'Arithmétique\/ $(\mathbb N,+,\times)$
avec paramètres\textup.
Alors\/ $\Th(M)$ est héréditairement indécidable\textup.
\end{proposition}
\begin{proof}
C'est un corollaire du théorème \ref{theorem:01.(6.3)} et des lemmes
\ref{lemma:07.(2.13)} et~\ref{lemma:19.(6.2)}.
\end{proof}

Le théorème C (voir l'introduction)
résulte maintenant du théorème A et
de la proposition \ref{proposition:04.(6.4)}.


\section{Questions ouvertes}
\label{section:questions}

Nous concluons avec deux questions qui, à notre connaissance,
sont ouvertes.

Le groupe $F$ de Thompson est définissable dans $T$.%
\footnote{Nous n'allons pas montrer cela ici.}
Par contre, il n'est pas connu aux auteurs si
$F$ est définissable dans~$V$.

\begin{question}
\label{question:1}
Le groupe $F$ de Thompson, est-il définissable dans $V$
avec paramètres ?
\end{question}


Nous avons déjà démontré que l'Arithmétique s'interprète dans $F$.
Par ailleurs, $F$ s'interprète dans l'Arithmétique parce que
le problème des mots y est décidable.
(Le théorème de Matiyasevich, voir
\cite{Matijasevich:1970:dpm-rus}
ou \cite[théorème 8.1]{Davis:1973:htpiu},
implique que tout sous-ensemble
récursif ou récursivement énumérable de $\mathbb N^n$, $n\in\mathbb N$,
est définissable dans l'Arithmétique ;
ainsi tout codage raisonnable des éléments de $F$
par des entiers positifs donne une interprétation de $F$
dans l'Arithmétique.)
Pourtant, même si une structure interprète l'Arithmétique et
que réciproquement elle est interprétée dans celle-ci,
il n'est pas nécessaire que les deux structures soient
\emph{bi-interprétables\/}
(voir \cite[théorème 6]{Khelif:2007:bisQFAegrac-fr} ou
\cite[théorème 7.16]{Nies:2007:dg}),
d'où la question.

\begin{question}
\label{question:2}
L'Arithmétique et le groupe $F$, sont-ils bi-interprétables parmi eux
avec paramètres ?
\end{question}

Disons que une structure $S$ est
\emph{catégoriquement finiment axiomatisée\/}
dans une classe $C$ de structures de la même signature
si $S\in C$ et qu'il existe un énoncé $\phi$ du premier ordre
tel que $S\models\phi$ et que toute structure dans $C$
qui satisfait $\phi$ soit isomorphe à $S$.%
\footnote{Dans le cas où la classe $C$ est formée de toutes
les structures de type fini d'une certaine classe,
André Nies \cite{Nies:2007:dg} et
Anatole Khélif \cite{Khelif:2007:bisQFAegrac-fr}
ont utilisé le terme 
«~\emph{quasi finiment axiomatisé\/}~»
dans plus ou moins le même sens que nous utilisons
«~catégoriquement finiment axiomatisé~».
Cependant les définitions dans \cite{Nies:2007:dg} et
\cite{Khelif:2007:bisQFAegrac-fr}
ne sont pas précises parce que la propriété d'être de type fini
n'est pas interne et dépend de la classe dans laquelle une structure
donnée est considérée.}
Selon Anatole Khélif \cite{Khelif:2007:bisQFAegrac-fr},
la bi-interprétabilité avec l'Arithmétique
peut être utilisée pour établir
l'axiomatisation finie catégorique dans des classes
de structures de signatures finie et de type fini.
Thomas Scanlon \cite{Scanlon:2008:ifgfbN}
a récemment établi la bi-interprétabilité
entre l'Arithmétique et tout corps commutatif de type fini
et a utilisé cela
pour démontrer que tout tel corps est
catégoriquement finiment axiomatisé dans la classe
de tous corps commutatifs de type fini, et ainsi
que la conjecture de Pop est vrai :
\begin{quote}
deux corps commutatifs de type fini sont élémentairement équivalents
si et seulement si ils sont isomorphes.
\end{quote}
André Nies a soulevé la question
s'il existe un groupe simple de type fini
catégoriquement finiment axiomatisé parmi
tous les groupes simples de type fini,
voir \cite[question 7.8]{Nies:2007:dg}.


\section{Appendice}
\label{section:appendix}

Ici nous démontrons que tout élément du sous-groupe dérivé de $\F$
est le produit de deux commutateurs, et donc que
$[\F,\F]$ est $0$-définissable dans $\F$.
La preuve de ce fait, que par ailleurs nous n'avons pas utilisé
dans ce travail, est connue aux spécialistes.
Cependant, ce fait est étroitement lié à la structure définissable
des groupes à l'étude,
et apparemment il ne paraît nulle part ailleurs
dans la littérature.
Pour le groupe $F$ de Thompson, ce résultat probablement fait
partie du folklore ; nous avons appris sa preuve de Matthew Brin,
qui avait légèrement modifié l'argument
de Keith Dennis et de Leonid Vaserstein
\cite[proposition 1(c)]{DennisVaserstein:1989:clg}.
Notre argument n'est qu'une généralisation triviale.

Comme dans la preuve de la proposition \ref{proposition:03.(5.2)},
notons
$$
\F^\circ=\{\,x\in\F\mid
(0)x^{\prime+}=(r)x^{\prime-}=1\,\}.
$$

\begin{proposition}
\label{proposition:05.(8.1)}
Tout élément de\/ $[\F,\F]$ est le produit de deux commutateurs
dans\/ $\F$\textup, voire dans\/ $\F^\circ$\/ \textup:
$$
[\F,\F]=\{\,[x_1,x_2][x_3,x_4]\mid x_1,x_2,x_3,x_4\in\F^\circ\,\}.
$$
\end{proposition}
\begin{proof}
Tout d'abord rappelons que $[\F,\F]\subset\F^\circ$
(voir le lemme \ref{lemma:10.(4.1)}).

Nous nous donnons maintenant deux éléments quelconques $x,y\in\F$
et leur commutateur $c=[x,y]$.
Soient $\alpha_1,\alpha_2,\beta_1,\beta_2\in A$ tels que
$0<\alpha_2<\alpha_1<\beta_1<\beta_2<r$ et que
$\supp(c)\subset]\alpha_1;\beta_1[$.
Alors il existe un endomorphisme $h\colon\F\to\F_{]\alpha_2;\beta_2[}$
tel que $h$ soit l'identité sur $\F_{]\alpha_1;\beta_1[}$.
En effet, soit $s\colon[0;r[\to[0;r[$ une application
qui est l'identité sur $[\alpha_1;\beta_1]$ et
qui est affine sur $[0;\alpha_1]$ et sur $[\beta_1;r[$
avec une pente $p\ll1$, $p\in\Lambda$,
alors on peut prendre comme $h$ la conjugaison par $s$
composée avec le plongement naturel des permutations de l'intervalle
$[(0)s;(r)s[$ dans les permutations de $[0;r[$
(cet $h$ est même injective).
Si $h$ est un tel endomorphisme, $(x)h=x'$, et que $(y)h=y'$, alors
$$
c=(c)h=[x',y']\quad\text{et}\quad
\supp(x')\cup\supp(y')\subset]\alpha_2;\beta_2[.
$$
Cela entraîne deux conséquences d'intérêt pour nous :
\begin{enumerate}
\item
    si $x,y\in\F$, alors il existe $x',y'\in\F^\circ$
    tels que $[x,y]=[x',y']$ ;
\item
    si $c_1,c_2,\dotsc,c_n$ sont des commutateurs dans $\F$,
    et $I_1,I_2,\dotsc,I_n$ sont des sous-intervalles
    deux-à-deux disjoints fermés de $[0;r[$, et
    tels que $\supp(c_i)\subset I_i$ pour tout $i=1,\dotsc,n$,
    alors le produit $c_1c_2\dotsm c_n$ est un commutateur aussi.
\end{enumerate}

Soient maintenant $c_1$, $c_2$, et $c_3$ trois commutateurs
arbitraires dans $\F$.
Choisissons $\alpha,\beta\in A\cap]0;r[$ tels que
$$
\supp(c_1)\cup\supp(c_2)\cup\supp(c_3)\subset]\alpha;\beta[.
$$
Choisissons $b\in\F$ tel que $(\alpha)b>\beta$.
Alors $0<(\beta)b^{-1}<\alpha<\beta<(\alpha)b<r$.
Comme
$$
\supp(c_2^b)\subset](\alpha)b;r[
\quad\text{et}\quad
\supp(c_3^{b^{-1}})\subset]0;(\beta)b^{-1}[,
$$
le produit $c_1c_2^bc_3^{b^{-1}}$ est un commutateur.%
\footnote{Nous utilisons la notation standard :
$x^y=y^{-1}xy$, $[x,y]=x^{-1}y^{-1}xy$.}
D'où
$$
c_1c_2c_3=c_1c_2^bc_3^{b^{-1}}
\Bigl((c_3^{-1})^{b^{-1}}c_2\Bigr)\Bigl((c_2^{-1})^bc_3\Bigr)
=(c_1c_2^bc_3^{b^{-1}})[c_2^{-1}c_3^{b^{-1}},b]
$$
est le produit de deux commutateurs.
\end{proof}




\newcommand{\noopsort}[1]{} \newcommand{\singleletter}[1]{#1}
  \providecommand{\href}[2]{#2} \providecommand{\url}[1]{\texttt{#1}}
  \providecommand{\nolinkurl}[1]{\texttt{#1}}
\providecommand{\bysame}{\leavevmode\hbox to3em{\hrulefill}\thinspace}
\providecommand{\MR}{\relax\ifhmode\unskip\space\fi MR }
\providecommand{\MRhref}[2]{%
  \href{http://www.ams.org/mathscinet-getitem?mr=#1}{#2}
}
\providecommand{\href}[2]{#2}
\begin{thebibliography}{10}

\bibitem{AhlbrandtZiegler:1986:qfatct}
Gisela Ahlbrandt et Martin Ziegler,
{\selectlanguage{english}%
  \emph{Quasi finitely axiomatizable totally categorical theories},
  Ann.\ Pure Appl.\ Logic \textbf{30} (1986), no.~1, 63--82}
  (anglais).

\bibitem{BardakovTolstykh:2007:iaTgF}
Valery~G.\ Bardakov et Vladimir~A.\ Tolstykh,
{\selectlanguage{english}%
  \emph{Interpreting the arithmetic in Thompson's group\/ $F$},
  J.\ Pure Appl.\ Algebra \textbf{211} (2007), no.~3, 633--637}
  (anglais).
  Prépublication :
   \href{http://arxiv.org/abs/math/0701748}{\nolinkurl{arXiv:math/0701748}}.

\bibitem{BelkBrown:2005:fdeTgF}
James~M.\ Belk et Kenneth~S.\ Brown,
{\selectlanguage{english}%
  \emph{Forest diagrams for elements of Thompson's group {$F$}},
  Internat.\ J.\ Algebra Comput.\ \textbf{15} (2005), no.~5--6, 815--850}
  (anglais).

\bibitem{BieriStrebel:pp1985:gPLhrl}
Robert Bieri et Ralph Strebel,
{\selectlanguage{english}%
  \emph{On groups of ${PL}$-homeomorphisms of the
  real line}},
  Math.\ Sem.\ der Univ.\ Frankfurt, notes non publiées, 1985.

\bibitem{BleakGGHMNS:pp2007:}
Collin Bleak, Alison Gordon, Garrett Graham, Jacob Hughes,
    Francesco Matucci, Hannah Newfield-Plunkett, et Eugenia Sapir,
{\selectlanguage{english}%
  \emph{Using dynamics to analyze centralizers in the generalized Higman-Thompson groups\/ $V_n$}},
  prépublication en version incomplète
  (anglais).

\bibitem{BrinSquier:2001:pcrcgplhrl}
Matthew~G.\ Brin et Craig.~C.\ Squier,
{\selectlanguage{english}%
  \emph{Presentations, conjugacy, roots, and centralizers in groups of
  piecewise linear homeomorphisms of the real line},
  Comm.\ Algebra, \textbf{29} (2001), no.~10, 4557--4596}
  (anglais).

\bibitem{Brown:1987:fpg}
Kenneth~S.\ Brown,
{\selectlanguage{english}%
  \emph{Finiteness properties of groups},
  J.\ Pure Appl.\ Algebra \textbf{44} (1987), no.~1--3, 45--75}
  (anglais).

\bibitem{CannonFloydParry:1996:inRTg}
James~W.\ Cannon, William~J.\ Floyd, et Walter~R.\ Parry,
{\selectlanguage{english}%
  \emph{Introductory notes on Richard Thompson's groups},
  Enseign.\ Math.\ (2) \textbf{42} (1996), no.~3--4, 215--256}
  (anglais).
  Prépublication :
    \url{http://www.geom.uiuc.edu/docs/preprints/lib/GCG63/thompson.ps}.

\bibitem{Davis:1973:htpiu}
Martin Davis,
{\selectlanguage{english}%
  \emph{Hilbert's tenth problem is unsolvable},
  Amer.\ Math.\ Monthly \textbf{80} (1973), 233--269}
  (anglais).

\bibitem{DelonSimonetta:1998:uwpspsf}
Fran{\c c}oise Delon et Patrick Simonetta,
{\selectlanguage{english}%
  \emph{Undecidable wreath products and skew power series fields},
  J.\ Symbolic Logic \textbf{63} (1998), no.~1, 237--246}
  (anglais).

\bibitem{DennisVaserstein:1989:clg}
R.~Keith Dennis et Leonid~N.\ Vaserstein,
{\selectlanguage{english}%
  \emph{Commutators in linear groups},
  $K$-Theory \textbf{2} (1989), no.~6, 761--767}
  (anglais).

\bibitem{Ershov:1980:prkm-rus}
Yuri~L.\ Ershov,
  \emph{Problemy razreshimosti i konstruktivnye modeli
    [Problèmes de décidabilité et modèles constructifs]},
    Matematicheskaya Logika i Osnovaniya Matematiki
    [Logique Mathématique et Fondations de Mathématique], 
  Nauka, Moscow, 1980
  (russe).


\bibitem{Godel:2006:fuSPMvS1-ger}
Kurt G{\"o}del,
{\selectlanguage{german}%
  \emph{{\"U}ber formal unentscheidbare {S}{\"a}tze der {P}rincipia
  {M}athematica und verwandter {S}ysteme. {I}
  [{S}ur les propositions formellement indécidables des
  Principia Mathematica et des systèmes apparentés. {I}]},
  Monatsh. Math. \textbf{149} (2006), no.~1, 1--30}
  (allemand),
  réimprimé sur
  Monatsh.\ Math.\ Phys.\ {\bf 38} (1931), 173--198,
  avec une introduction par Sy-David Friedman.

\bibitem{Higman:1974:fpisg}
Graham Higman,
{\selectlanguage{english}%
  \emph{Finitely presented infinite simple groups},
  Notes on Pure Mathematics, vol.~8,
  Australian National University, Canberra, 1974}
  (anglais).

\bibitem{Hodges:1993:mt}
Wilfrid Hodges,
{\selectlanguage{english}%
  \emph{Model theory},
  Encyclopedia of Mathematics and its Applications, vol.~42,
  Cambridge University Press, 1993}
  (anglais).

\bibitem{Hodges:1997:smt}
\bysame,
{\selectlanguage{english}%
  \emph{A shorter model theory},
  Cambridge University Press, 1997}
  (anglais).

\bibitem{Khelif:2007:bisQFAegrac-fr}
Anatole Kh{\'e}lif,
  \emph{Bi-interpr{\'e}tabilit{\'e}e et structures {QFA} :
    {\'e}tude de groupes r{\'e}solubles et des anneaux commutatifs},
  C.\ R.\ Math.\ Acad.\ Sci.\ Paris \textbf{345} (2007),
  no.~2, 59--61.

\bibitem{Matijasevich:1970:dpm-rus}
Yuri~V.\ Matijasevich,
  \emph{Diofantovost' perechislimykh mnozhestv
    [La nature diophantienne des ensembles énumérables]},
  Dokl.\ Akad.\ Nauk SSSR \textbf{191} (1970), 279--282
  (russe).
  Traduction anglaise dans Soviet Math.\ Dokl.

\bibitem{MostowskiRobTar:1971:ueua}
Andrzej Mostowski, Raphael~M.\ Robinson, et Alfred Tarski,
{\selectlanguage{english}%
  \emph{Undecidability and essential undecidability in arithmetic},
  Undecidable theories,
  Studies in logic and the foundations of mathematics,
  North-Holland Publishing Co., Amsterdam, 1971, pp.~36--74}
  (anglais).

\bibitem{MostowskiTarski:1949:uaitr}
Andrzej Mostowski et Alfred Tarski,
{\selectlanguage{english}%
  \emph{Undecidability in the arithmetic of
    integers and in the theory of rings},
  J.\ Symbolic Logic \textbf{14} (1949), 76}
  (anglais).

\bibitem{Nies:2007:dg}
Andr{\'e} Nies,
{\selectlanguage{english}%
  \emph{Describing groups},
  Bull.\ Symbolic Logic \textbf{13} (2007), no.~3, 305--339}
  (anglais).

\bibitem{Noskov:1983:etkpprg-rus}
Gennady~A.\ Noskov,
  \emph{Ob elementarnoj teorii konechno porozhdyennoj pochti
  razreshimoj gruppy
  [À propos de la théorie élémentaire d'un groupe virtuellement
  résoluble de type fini]},
  Izv.\ Akad.\ Nauk SSSR Ser.\ Mat.\ \textbf{47} (1983), no.~3, 498--517
  (russe).
  Traduction anglaise dans Math.\ USSR Izvestiya.

\bibitem{Noskov:1984:etfgasg-eng}
\bysame,
{\selectlanguage{english}%
  \emph{On the elementary theory of a finitely generated
    almost solvable group},
  Math.\ USSR Izvestiya \textbf{22} (1984), no.~3, 465--482}
  (anglais).
  Traduit du russe.

\bibitem{Pillay:1983:ist}
Anand Pillay,
{\selectlanguage{english}%
  \emph{An introduction to stability theory},
  Oxford Logic Guides, vol.~8,
  The Clarendon Press, Oxford University Press, New York, 1983}
  (anglais).

\bibitem{Poizat:1985:ctm-fr}
Bruno Poizat,
  \emph{Cours de th{\'e}orie des mod{\`e}les.
    {U}ne introduction {\`a} la logique math{\'e}matique contemporaine},
  Nur al-Mantiq wal-Ma'rifah, Bruno Poizat, Lyon, 1985.

\bibitem{Poizat:1987:gs-fr}
\bysame, 
  \emph{Groupes stables.
    Une tentative de conciliation entre la g{\'e}om{\'e}trie alg{\'e}brique
    et la logique math{\'e}matique},
  Nur al-Mantiq wal-Ma'rifah, vol.~2, Bruno Poizat, Lyon, 1987.

\bibitem{Poizat:2001:sg-eng}
\bysame,
{\selectlanguage{english}%
  \emph{Stable groups},
  Mathematical Surveys and Monographs, vol.~87,
  American Mathematical Society, 2001}
  (anglais).
  Traduit de l'original français de 1987 par Moses Gabriel Klein.

\bibitem{Prest:1988:mtm}
Mike~Y.\ Prest,
{\selectlanguage{english}%
  \emph{Model theory and modules},
  London Mathematical Society
  Lecture Note Series, vol.\ 130,
  Cambridge University Press, 1988}
  (anglais).

\bibitem{Robinson:1951:ur}
Raphael~M.\ Robinson,
{\selectlanguage{english}%
  \emph{Undecidable rings},
  Trans.\ Amer.\ Math.\ Soc.\ \textbf{70} (1951), 137--159}
  (anglais).

\bibitem{Rothmaler:2000:imt-eng}
Philipp Rothmaler,
{\selectlanguage{english}%
  \emph{Introduction to model theory},
  Algebra, Logic and Applications, vol.~15,
  Gordon and Breach Science Publishers, Amsterdam, 2000}
  (anglais).
  Traduit et révisé de l'original allemand de 1995 par l'auteur.

\bibitem{Scanlon:2008:ifgfbN}
Thomas Scanlon,
{\selectlanguage{english}%
  \emph{Infinite finitely generated fields are
    biinterpretable with\/ $\mathbb {N}$}},
  J.\ Amer.\ Math.\ Soc. \textbf{21} (2008), no.~3, 893--908,
  (anglais).

\bibitem{Sela:pp2006:dgg8s}
Zlil Sela,
{\selectlanguage{english}%
  \emph{Diophantine geometry over groups VIII: stability}},
  prépublication,
  \href{http://arxiv.org/abs/math/0609096}{\nolinkurl{arXiv:math/0609096}},
  2006
  (anglais).

\bibitem{Stein:1992:gplh}
Melanie Stein,
{\selectlanguage{english}%
  \emph{Groups of piecewise linear homeomorphisms},
  Trans.\ Amer.\ Math.\ Soc.\ \textbf{332} (1992), no.~2, 477--514}
  (anglais).

\bibitem{Tarski:1971:gmpu}
Alfred Tarski,
{\selectlanguage{english}%
  \emph{A general method in proofs of undecidability},
  Undecidable theories,
  Studies in logic and the foundations of mathematics,
  North-Holland Publishing Co., Amsterdam, 1971, pp.~1--35}
  (anglais).

\bibitem{Wagner:2000:sg}
Frank~O.\ Wagner,
{\selectlanguage{english}%
  \emph{Stable groups},
  Handbook of algebra, vol.~2,
  North-Holland Publishing Co., Amsterdam, 2000}
  (anglais).

\bibitem{anonymous:2004:Tg40y}
\emph{Thompson's Group at 40 Years},
{\selectlanguage{english}%
  Problem list of the workshop held 11--14 Jan., 2004,
  at American Institute of Mathematics, Palo Alto, California},
  \url{http://www.aimath.org/WWN/thompsonsgroup/thompsonsgroup.pdf}, 2004
  (anglais).

\end{thebibliography}
%

\newcommand{\noopsort}[1]{}
\newcommand{\singleletter}[1]{#1}
\providecommand{\href}[2]{#2}
\providecommand{\url}[1]{\texttt{#1}}
\providecommand{\nolinkurl}[1]{\texttt{#1}}
\providecommand{\bysame}{\leavevmode\hbox to3em{\hrulefill}\thinspace}
\providecommand{\MR}{\relax\ifhmode\unskip\space\fi MR }
\providecommand{\MRhref}[2]{%
  \href{http://www.ams.org/mathscinet-getitem?mr=#1}{#2}
}



\end{document}